%% file: HD_Sobolev_arXiv_v2.tex
\documentclass[oneside,11pt]{article}

\input{preamble_arXiv.tex}

\usepackage[pdftex,final,bookmarksnumbered,bookmarksopen=false,breaklinks,colorlinks]{hyperref}
\hypersetup{
	final,
	bookmarksnumbered=true,
	bookmarksopen=true,
	bookmarksopenlevel=0,
	unicode=false,
	pdftoolbar=true,
	pdfmenubar=true,
	pdffitwindow=false,
	pdftitle={High-dimensional Sobolev tests on hyperspheres},
	pdfdisplaydoctitle=true,
	pdftoolbar=true,
	pdfmenubar=true,
	pdflang={English},
	pdfauthor={Bruno Ebner, Eduardo Garcia-Portugues, Thomas Verdebout},
	pdfsubject={arXiv paper},
	pdfcreator={Eduardo Garcia-Portugues},
	pdfproducer={Eduardo Garcia-Portugues},
	pdfkeywords={Directional statistics; Normality; Rotational symmetry; Spherical symmetry; Uniformity},
	pdfnewwindow=true,
	breaklinks=true,
	hidelinks,
	linkcolor=black,
	citecolor=black,
	filecolor=magenta,
	urlcolor=cyan
}

\newif\ifmain
\maintrue
\newif\ifsupplement
\supplementtrue

\newif\iffigstabs
\figstabstrue

\begin{document}

\ifmain

\title{High-dimensional Sobolev tests on hyperspheres}
\setlength{\droptitle}{-1cm}
\predate{}%
\postdate{}%
\date{}

\author{Bruno Ebner$^{1}$, Eduardo Garc\'ia-Portugu\'es$^{2,4}$, and Thomas Verdebout$^{3}$}
\footnotetext[1]{Institute of Stochastics, Karlsruhe Institute of Technology (Germany).}
\footnotetext[2]{Department of Statistics, Universidad Carlos III de Madrid (Spain).}
\footnotetext[3]{Department of Mathematics, Université libre de Bruxelles (Belgium).}
\footnotetext[4]{Corresponding author. e-mail: \href{mailto:edgarcia@est-econ.uc3m.es}{edgarcia@est-econ.uc3m.es}.}
\maketitle

\begin{abstract}
We derive the limit null distribution of the class of Sobolev tests of uniformity on the hypersphere as the dimension and the sample size diverge to infinity at arbitrary rates. The limiting non-null behavior of these tests is also explored: (\textit{i}) a general consistency result for Sobolev tests in high dimensions is shown and (\textit{ii}) the asymptotic power of the tests under sequences of integrated von Mises--Fisher local alternatives is obtained. The asymptotic results are applied to test for high-dimensional rotational symmetry and spherical symmetry. Numerical experiments illustrate the derived behavior of the uniformity and spherical symmetry tests under the null and under local and fixed alternatives. A real data application tests the high-dimensional normality of the cosmic microwave background.
\end{abstract}
\begin{flushleft}
\end{flushleft}

\section{Introduction}
\label{sec:intro}

Testing the uniformity of a sample of $n$ independent and identically distributed random vectors $\bX_{1}, \ldots, \bX_{n}$ on the unit (hyper)sphere $\Sd\defin\{\bx\in\R^{d+1}:\bx^\top\bx=1\}$, $d\geq1$, is a classical testing problem in statistics and machine learning. Detecting whether a sample departs from spherical uniformity is indeed fundamental in problems ranging from goodness-of-fit assessment to the analysis of directional data, random projections, and geometric properties of learned representations. As a result, the problem has generated many contributions in the statistical literature: a number of available tests are described in \cite{Garcia-Portugues2020a}. High-dimensional directional data have been considered in~\cite{Dryden2005}, with applications in brain shape modeling, and in \cite{Banerjee2005}, with applications to text mining, to name but a few. Within uniformity testing, \cite{Cuesta-Albertos2009} proposed a projection-based test that performs well empirically even in high dimensions, \cite{Paindaveine2016} studied the high-dimensional null behavior of the Rayleigh test, while \cite{Cutting2017, Cutting2021} focused on the high-dimensional asymptotic power of the Rayleigh and Bingham tests. \cite{Chikuse1991a,Chikuse1993} also explicitly considered high-dimensional testing for uniformity on the sphere in a fixed-$n$ large-$d$ framework, while \cite{Cai2013} rather adopted a double asymptotic approach for the same problem. In machine learning, \cite{Liu2018learning} demonstrated that neural networks trained to minimize hyperspherical energy naturally produce standardized neurons whose distribution converges to the uniform distribution on the hypersphere. Therefore, for overparametrized neural networks, testing for uniformity on high-dimensional spheres can offer valuable insights into the structural regularization of layers.

One of the most overarching frameworks for testing uniformity on $\Sd$ is the class of \emph{Sobolev tests}, a term coined by \cite{Gine1975} to reflect their relation with Sobolev norms. This class of tests was introduced in \cite{Beran1968} and \cite{Gine1975}, and contains the famous \cite{Rayleigh1919} and \cite{Bingham1974} tests, respectively based on the sample average and sample covariance matrix. Sobolev tests reject the null hypothesis of uniformity for large values of statistics of the form
\begin{align} \label{Sobol}
    S_n\defin \frac{1}{n} \sum_{i,j=1}^n\sum_{k=1}^\infty v_{k,d}^2 h_k(\bX_{i}^\top \bX_{j}),\quad h_k(x)\defin \begin{cases}
    2 \cos\big(k \cos^{-1}(x)\big), & \text{if } d=1, \\
    \bigg(1+ \frac{2k}{d-1}\bigg) C_k^{(d-1)/2} (x), & \text{if } d>1,
    \end{cases}
\end{align}
where $(v_{k,d})_{k=1}^\infty$ is a real sequence and $C_k^{(d-1)/2}$ is a Gegenbauer polynomial (see Section \ref{sec:asymp}). Data-driven Sobolev tests have been proposed in \cite{Bogdan2002} for the circular case, and on compact Riemannian manifolds in \cite{Jupp2008, Jupp2009}. The low-dimensional asymptotic behavior of Sobolev tests has been studied in \cite{Gine1975}, while detection thresholds under rotationally symmetric alternatives have been obtained recently in \cite{Garcia-Portugues:Sobolev}.

The high-dimensional limiting behavior of (general) Sobolev tests remains largely unexplored; \cite{Paindaveine2016, Cutting2017, Cutting2021} and \cite{jiang2025asymptotic, jiang2025detecting} considered the asymptotic behavior of specific Sobolev tests (the Rayleigh test, the Bingham test, etc.) in various settings but did not consider the entire class of Sobolev tests. The present paper fills this gap by deriving the asymptotic null distributions of a general Sobolev statistic on $\Sdn$ as the dimension $d_n$ diverges to infinity together with the sample size $n$. Additionally, the high-dimensional consistency against general alternatives and the high-dimensional non-null behavior of Sobolev tests against integrated von Mises--Fisher local alternatives are derived. We highlight that no assumptions on the rates of growth of $n$ and $d_n$ are imposed in these two results, and that they generalize the null and non-null high-dimensional findings by \cite{Paindaveine2016} and \cite{Cutting2017} for the Rayleigh test. Applications of the high-dimensional limit results are then presented for testing rotational symmetry on the hypersphere $\Sdn$ and testing spherical symmetry on $\R^{d_n+1}$. Numerical experiments in the form of Monte Carlo simulations are conducted to corroborate the asymptotic findings, and the methodology is illustrated on a real data application testing the high-dimensional normality of the cosmic microwave background. The appendices contain the proofs, and the Supplementary Materials (SM) provide complementary numerical experiments.

The rest of this paper is organized as follows. Section \ref{sec:asymp} contains the main results on the null and non-null asymptotic behavior of high-dimensional Sobolev tests of uniformity on $\Sdn$. Applications to testing high-dimensional rotational symmetry on $\Sdn$ and spherical symmetry on $\R^{d_n+1}$ are given in Section \ref{sec:ext}. Section \ref{sec:num} provides a simulation study on testing high-dimensional spherical symmetry, and Section \ref{sec:cmb} illustrates the tests on the cosmic microwave background. A discussion is given in Section \ref{sec:disc}. Appendices \ref{app:proofs} and \ref{app:lemmas} collect the proofs of the presented results, while Appendix C
in the SM shows an empirical illustration of the main asymptotic results.

The SM and the repository \url{https://github.com/egarpor/sphunif/} contain the R scripts to reproduce the simulations, numerical experiments, and real data application.

\section{Main results}
\label{sec:asymp}

\subsection{General setup}
\label{sec:general}

Let $\bX_{n,i}$, $i=1,\ldots,n$, $n=1,2,\ldots$, be a triangular array of random vectors such that, for any~$n$, the~$\bX_{n,i}$'s are mutually independent and share a common distribution on the unit sphere~$\Sdn$. We consider the problem of testing the null hypothesis~$\mathcal{H}_{0,n}$ stating that this common distribution is the uniform distribution on~$\Sdn$ in the asymptotic regime where both $n$ and $d_n$ diverge to infinity.

The Sobolev statistics are based on the Gegenbauer polynomials $\{C_k^{(d-1)/2}\}_{k=0}^\infty$ of index $(d-1)/2$, $d\geq2$, which form an orthogonal basis of the space $L_d^2([-1,1])$ of square-integrable real functions on $[-1,1]$ with respect to the weight $x\mapsto (1-x^2)^{d/2-1}$. More precisely, they satisfy the orthogonality relation
\begin{align}
	\int_{-1}^1 C_k^{(d-1)/2}(x)C_\ell^{(d-1)/2}(x)(1-x^2)^{d/2-1}\,\rd x=\delta_{k,\ell}c_{k,d},\label{eq:ckd1}
\end{align}
where $\delta_{k,\ell}$ denotes the Kronecker delta and
\begin{align}
    c_{k,d}&\defin\frac{2^{3-d}\pi\Gamma(d+k-1)}{(d+2k-1)k!\Gamma((d-1)/2)^2}=\frac{\om{d}}{\om{d-1}}\lrp{1+\frac{2k}{d-1}}^{-2}d_{k,d},\label{eq:ckd2}
\end{align}
with $\om{d}\defin 2 \pi^{(d+1)/2}/ \Gamma((d+1)/2)$ denoting the surface area of $\Sd$. Above,
\begin{align}
    d_{k,d}&\defin \binom{d+k-2}{d-1}+\binom{d+k-1}{d-1}=\lrp{1+\frac{2k}{d-1}}\frac{\Gamma(d-1+k)}{\Gamma(d-1)k!}\label{eq:dkd}
\end{align}
is the dimension of the vector space that contains the collection of spherical harmonics of degree $k$, i.e., the space of harmonic homogeneous polynomials of degree $k$ defined on $\Sd$. As $k\to\infty$, $d_{k,d}\asymp k^{d-1}$, with $a_k\asymp b_k$ denoting $a_k/b_k\to c\neq0$ as $k\to\infty$. The first Gegenbauer polynomials are $C_0^{(d-1)/2}(x)=1$, $C_1^{(d-1)/2}(x)=(d-1)x$, and $C_2^{(d-1)/2}(x)=((d^2-1)x^2 -d+1)/2$, while the subsequent ones follow from the recurrence relation
\begin{align*}
    (k+1)C_{k+1}^{(d-1)/2}(x)=(2k+d-1)\,x\,C_k^{(d-1)/2}(x)-(k+d-2)C_{k-1}^{(d-1)/2}(x).
\end{align*}
It also follows that $C_k^{(d-1)/2}(1)=\lrp{1+2k/(d-1)}^{-1}d_{k,d}$.

When $d=1$, Sobolev statistics are based on Chebyshev polynomials (of the first kind) $\{C_k^{0}\}_{k=0}^\infty$, defined by $C_k^0(x)\defin \cos(k\cos^{-1}(x))$. They can be seen as the scaled limit of Gegenbauer polynomials when their index shrinks to zero: $\lim_{\lambda\to0^+}\lambda^{-1}C_k^{\lambda}(x)=(2/k)C_k^{0}(x)$. Chebyshev polynomials form an orthogonal basis in $L_1^2([-1,1])$, the space of square-integrable functions in $[-1,1]$ with respect to $x\mapsto (1-x^2)^{-1/2}$, satisfying
\begin{align*}
    \int_{-1}^1 C_k^{0}(x)C_\ell^{0}(x)(1-x^2)^{-1/2}\,\rd x=\delta_{k,\ell}c_{k,1},\quad c_{k,1}\defin\frac{1+\delta_{k,0}}{2}\pi.
\end{align*}
In this case, $d_{k,1}=2 - \delta_{k,0}$ is the dimension of the space of harmonic homogeneous polynomials of degree $k$ on $\mathbb{S}^1$.

Any function $g\in L_d^2([-1,1])$, $d\geq 1$, admits the representation
\begin{align*}
    g(x)=\sum_{k=0}^{\infty} b_{k,d}(g) \,C_{k}^{(d-1)/2}(x),\quad b_{k,d}(g)\defin \frac{1}{c_{k,d}}\int_{-1}^1 g(x)C_k^{(d-1)/2}(x)(1-x^2)^{d/2-1}\,\rd x
\end{align*}
with the series converging in mean square.

\subsection{Null behavior}

In our first main result, we study the high-dimensional asymptotic behavior under $\mathcal{H}_{0,n}$ of a general Sobolev test statistic
$$
T_n\defin\frac{2}{n}\sum_{1\leq i<j\leq n} \psin(\bX_{n,i}^\top \bX_{n,j})
$$
with $\psin:[-1,1]\to\mathbb{R}$ defined almost everywhere as
\begin{align}
	\psin(x)\defin\sum_{k=1}^\infty \lrp{1+\frac{2k}{d_n-1}}v_{k,d_n}^2 C_k^{(d_n-1)/2}(x),\label{eq:psi}
\end{align}
for a real sequence of weights $(v_{k,d_n})_{k=1}^\infty$ and $d_n\geq2$. For $d_n=1$, let $$
\psin(x)\defin\sum_{k=1}^\infty 2v_{k,1}^2 C_k^0(x).
$$
Observe that the test statistic $T_n$ is the $U$-statistic form of \eqref{Sobol}.

\begin{theorem} \label{thm:asymp}
Let $d_n$ be a sequence of positive integers such that $d_n\to\infty$ as $n\to\infty$. Assume that $\bX_{n,i}$, $i=1,\ldots,n$, $n=1,2,\ldots$, form a triangular array such that, for any fixed $n$, $\bX_{n,1},\ldots, \bX_{n,n}$ are mutually independent and uniformly distributed on $\Sdn$. Assume moreover that:
\begin{enumerate}[label=(\textit{\roman{*}}), ref=(\textit{\roman{*}})]
    \item For every $n\geq 1$ and $d_n\geq 1$, \label{cond1}
    \begin{align}
        \sigma_{n}^2\defin 2\E{\psin^2(\bX_{n,1}^\top\bX_{n,2})}=2\sum_{k=1}^\infty v_{k,d_n}^4 d_{k,d_n}<\infty.\label{eq:sigman2}
    \end{align}
    \item As $n\to\infty$ and $d_n\to\infty$, \label{cond2}
    \begin{align}
        \sigma_{n}^{-4}\sum_{k=1}^\infty v^8_{k,d_n} d_{k,d_n}=o(1)\quad\text{and}\quad \sigma_{n}^{-4}\Ebig{\psin^4(\bX_{n,1}^\top\bX_{n,2})}=O(1). \label{eq:Av8}
    \end{align}
\end{enumerate}

Then $\sigma_{n}^{-1} T_n\inlaw\mathcal{N}(0,1)$, with $\inlaw$ denoting weak convergence.
\end{theorem}

\begin{remark} \label{rem:finvar}
Condition \ref{cond1} refers to the finiteness of the asymptotic variance of $T_n$, and implies that, when restricted to power orders, the sequence of weights is such that $v_{k,d}^2=O(k^{-(d+\delta)/2})$ as $k\to\infty$ for $\delta>0$. This requirement is typically met by most Sobolev test statistics, but there are some exceptions. For example, the stereographic test \citep{Fernandez-de-Marcos:stereo} has a weight sequence $v_{k,d}^2\asymp k^{-d+1}$, which makes $\sigma_{n}^2$ finite only if $d_n\geq 3$. For this and other statistics where $\sigma_n^2<\infty$ only if $d_n\geq d_*$, Theorem \ref{thm:asymp} is still applicable to the new triangular array $\tilde{\bX}_{n,i}\defin\bX_{n+n_*-1,i+n_*-1}$, with $n_*\in\mathbb{N}$ such that $d_{n}\geq d_*$ for $n\geq n_*$.
\end{remark}

\begin{remark} \label{rem:psi4}
The first part of condition \ref{cond2} is implied by the simpler condition
\begin{align*}
    \sigma_n^{-2}\max_{k\geq 1} v_{k,d_n}^4 = o(1)
\end{align*}
as $n\to\infty$, since $\sum_{k=1}^\infty v_{k,d_n}^8 d_{k,d_n}
\leq \big(\max_{k\geq 1} v_{k,d_n}^4\big)\sum_{k=1}^\infty v_{k,d_n}^4 d_{k,d_n}
= \frac{1}{2}\big(\max_{k\geq 1} v_{k,d_n}^4\big)\,\sigma_n^2$. The ratio in the second part of \ref{cond2} satisfies
\begin{align*}
    \sigma_n^{-4}\Ebig{\psin^4(\bX_{n,1}^\top\bX_{n,2})}\geq \frac{1}{4},
\end{align*}
for all $n,d_n\geq1$, due to Jensen's inequality: $\Ebig{\psin^4(\bX_{n,1}^\top\bX_{n,2})}\geq \Ebig{\psin^2(\bX_{n,1}^\top\bX_{n,2})}^2=\sigma_n^4/4$. As such, the $O(1)$ in \ref{cond2} cannot be $o(1)$.
\end{remark}

It is interesting to apply Theorem \ref{thm:asymp} to the case of a ``$k_0$-Sobolev test'' where, for a fixed integer $k_0\geq1$, $v_{k,d_n}=\delta_{k,k_0}$, for $k\geq1$. Then,
$$
\psi_{n}(x)=\psi_{n,k_0}(x)\defin\left(1+\frac{2k_0}{d_n-1}\right)C_{k_0}^{(d_n-1)/2}(x)
$$
and, because of Lemma \ref{lem:Ckm},
$\Ebig{\psi_{n,k_0}^4(\bX_{n,1}^\top\bX_{n,2})}
 \asymp d_n^{2k_0}$.
Therefore, conditions \ref{cond1} and \ref{cond2} hold, since $\sigma_n^2=2d_{k_0,d_n}\asymp d_n^{k_0}<\infty$
for every $n\geq 1$ and $d_n\geq 1$, and
$$
\sigma_n^{-4}\sum_{k=1}^\infty \delta_{k,k_0} d_{k,d_n}\asymp d_n^{-2k_0}d_n^{k_0}=o(1)\quad\text{and}\quad \sigma_n^{-4}\Ebig{\psin^4(\bX_{n,1}^\top\bX_{n,2})}\asymp d_n^{-2k_0}d_n^{2k_0}=O(1).
$$
Therefore, Theorem \ref{thm:asymp} gives
\begin{align}
	\sqrt{\frac{2}{d_{k_0,d_n}}}\left(1+\frac{2k_0}{d_n-1}\right) \frac{1}{n}\sum_{1\leq i<j\leq n} C_{k_0}^{(d_n-1)/2}(\bX_{n,i}^\top\bX_{n,j})\inlaw \mathcal{N}(0,1)\label{eq:asympCkexact}
\end{align}
as $n\to\infty$ and $d_n\to\infty$. Since $d_{k_0,d_n}\sim\lrp{1+2k_0/(d_n-1)}d_n^{k_0}/k_0!$ by Lemma \ref{lem:Ckm} ($\sim$ denoting equality in the limit), a simpler asymptotic equivalence of \eqref{eq:asympCkexact} is
\begin{align}
	\frac{\sqrt{2\,k_0!}}{d_n^{k_0/2}n}\sum_{1\leq i<j\leq n} C_{k_0}^{(d_n-1)/2}(\bX_{n,i}^\top\bX_{n,j})\inlaw \mathcal{N}(0,1)\label{eq:asympCk}
\end{align}
as $n\to\infty$ and $d_n\to\infty$. For $k_0=1$, the high-dimensional result in \citet[Theorem 2.1]{Paindaveine2016} for the Rayleigh statistic arises as a particular case in \eqref{eq:asympCk}:
\begin{align*}
    \frac{\sqrt{2d_n}}{n}\sum_{1\leq i<j\leq n} \bX_{n,i}^\top\bX_{n,j}
    \sim
    \frac{\sqrt{2}}{d_n^{1/2}n}\sum_{1\leq i<j\leq n} C_{1}^{(d_n-1)/2}(\bX_{n,i}^\top\bX_{n,j})
    \inlaw \mathcal{N}(0,1).
\end{align*}

\subsection{Non-null behavior: consistency}

In the present section, we discuss some high-dimensional power properties of our tests. First, note that, as explained in \citet[Section 3]{Garcia-Portugues:Sobolev}, we have that for $\bu,\bv \in \Sdn$,
\begin{align} \label{Harmo}
    \lrp{1+\frac{2k}{d_n-1}}C_k^{(d_n-1)/2}(\bu^\top\bv)= \sum_{r=1}^{d_{k,d_n}} g_{r, k}(\bu) g_{r,k}(\bv),
\end{align}
where $g_{1,k}, \ldots, g_{d_{k, d_n},k}$ form an orthonormal basis of the linear subspace of real homogeneous polynomials of degree $k$ whose Laplacian is zero. Let $\bX_{n,i}$, $i=1,\ldots,n$, $n=1,2,\ldots$, form a triangular array such that, for any fixed $n$, $\bX_{n,1},\ldots, \bX_{n,n}$ are mutually independent and identically distributed with common distribution $F_n$. We have the following result on the test $\phi_n$ that rejects the null hypothesis at the nominal level $\alpha \in (0,1)$ when $\sigma_n^{-1} T_n > z_{\alpha}$, where $z_{\alpha}$ is the upper $\alpha$-quantile of the standard Gaussian distribution.

\begin{theorem}\label{thm:consis}
Assume that $\bX_{n,1},\ldots, \bX_{n,n}$ are mutually independent and identically distributed with common distribution $F_n$ such that $e_{r,k,n}\defin\mathrm{E}_{F_n}[g_{r, k}(\bX_{n,1})]$ and
$$
c_{r,r',n}^{(k,k')}\defin \mathrm{Cov}_{F_n}(g_{r, k}(\bX_{n,1}), g_{r',k'}(\bX_{n,1}))
$$
are well-defined and satisfy
\begin{enumerate}[label=(\textit{\roman{*}}), ref=(\textit{\roman{*}})]
    \item $\sigma_n^{-1}\sum_{k=1}^\infty v_{k,d_n}^2 \sum_{r=1}^{d_{k,d_n}} e_{r,k,n}^2 \to c>0$,\label{thm:consis:1}
    \item $\sum_{k, k'=1}^\infty v_{k,d_n}^2 v_{k',d_n}^2 \sum_{r=1}^{d_{k,d_n}} \sum_{r'=1}^{d_{k',d_n}}  e_{r,k,n} e_{r',k',n} c_{r,r',n}^{(k,k')}= o(\sigma_n^{2} n)$, and\label{thm:consis:2}
    \item $\sum_{k, k'=1}^\infty v_{k,d_n}^2 v_{k',d_n}^2 \sum_{r=1}^{d_{k,d_n}} \sum_{r'=1}^{d_{k',d_n}} \big(c_{r,r',n}^{(k,k')}\big)^2= o(\sigma_n^{2} n^{2})$ \label{thm:consis:3}
\end{enumerate}
as $n \to \infty$. Then, we have that $$\lim_{n \to \infty} \mathrm{E}_{F_n}[\phi_{n}]=1.$$
\end{theorem}

We now discuss some consequences of the result. Although the result is very general, we focus here on rotationally symmetric alternatives in order to situate the above result within the existing body of work in this setting. As explained in Section \ref{sec:sphsym} below, a random vector~$\bX$ on~$\Sdn$ has a \emph{rotationally symmetric} distribution about $\btheta \in \Sdn$ if and only if~$\bO \bX$ has the same distribution as $\bX$ for any orthogonal matrix~$\bO\in\mathrm{SO}(d_n+1)$ such that $\bO\btheta=\btheta$. Note that when we consider the asymptotic power of Sobolev tests, we can safely fix any value of $\btheta$ thanks to the rotation-invariance property of Sobolev tests; we therefore fix $\btheta=(1, 0, \ldots, 0)^\top \in \Sdn$. Under rotational symmetry about $\btheta$, we have that $e_{r,k,n}=0$ for $r>1$, with $g_{1,k}$ the zonal harmonic about $\btheta$ (see \citet{Garcia-Portugues:Sobolev}). For a $k_0$-Sobolev test where, for a fixed integer $k_0\geq1$, $v_{k,d_n}=\delta_{k,k_0}$, for $k\geq1$, conditions \ref{thm:consis:1}, \ref{thm:consis:2}, and \ref{thm:consis:3} are respectively equivalent (up to unimportant constants) to
\begin{enumerate}[label=(\textit{\roman{*}}), ref=(\textit{\roman{*}})]
    \item $d_{k_0,d_n}^{-1/2} e_{1,k_0,n}^2 \to c$, \label{conddk:1}
    \item $d_{k_0,d_n}^{-1} n^{-1}e_{1,k_0,n}^2 c_{1,1,n}^{(k_0,k_0)} \to 0$ as $n \to \infty$, and \label{conddk:2}
    \item $d_{k_0,d_n}^{-1} n^{-2} \sum_{r,r'=1}^{d_{k_0,d_n}}\big(c_{r,r',n}^{(k_0,k_0)}\big)^2 \to 0$ as $n \to \infty$. \label{conddk:3}
\end{enumerate}
In the low-dimensional case ($d_n$ fixed), these conditions really read as classical fixed alternatives. Indeed, taking for instance $k_0=1$ (the Rayleigh test) and the von Mises--Fisher alternative with density $\bx \mapsto c_{d_n,\kappa_n} \exp(\kappa_n\, \bx^\top\btheta)$, where $\kappa_n \geq 0$ and $\btheta \in \Sdn$, condition \ref{conddk:1} is equivalent to non-vanishing values of $e_{1,1,n}$ that can be obtained for non-vanishing values of $\kappa_n$ (and therefore fixed alternatives).

Obviously, the high-dimensional version of the problem is more complicated from the statistical point of view, and the conditions above can be seen as ``high-dimensional'' fixed alternatives under which $k_0$-Sobolev tests are consistent. Since $d_{k_0,d_n}\sim d_n^{k_0}/k_0!$, we obtain for the Rayleigh test ($k_0=1$) that if $e_{1,1,n} \sim d_n^{1/4}$ and $c_{1,1,n}^{(1,1)}=o(nd_n^{1/2})$ as $n \to \infty$, the Rayleigh test is consistent as $n\to\infty$ and $d_n \to \infty$. Here \ref{conddk:3} reduces to the latter condition because, for $k_0=1$, the non-zonal harmonic variances are bounded and hence contribute only $O(n^{-2})$. Therefore, the consistency of the Rayleigh test is obtained for sequences of rotationally symmetric distributions such that $e_{1,1,n} \sim d_n^{1/4}$ and therefore $e_{1,1,n} \to \infty$ as $n \to \infty$. This reflects the fact that the high-dimensional problem is statistically more challenging than the low-dimensional one, which is not particularly surprising. In the next section, we move beyond consistency results. In the von Mises--Fisher case, we provide more informative results on the detection threshold of Sobolev tests.

\subsection{Non-null behavior: the integrated von Mises--Fisher case}

We now explore the situation where, for any $n$, the sample $\bX_{n,1},\ldots,\bX_{n,n}$ follows an \emph{integrated von Mises--Fisher distribution} on $\Sdn$ with location parameter $\btheta\in\Sdn$ and varying concentration parameter $\kappa_n\geq 0$ (with $\kappa_n=0$ corresponding to uniformity). Such alternatives are characterized by likelihood functions (with respect to the product surface area measure on~$(\Sdn)^n$) of the form
\begin{align}
    \frac{\rd{\rm P}^{(n)}_{\kappa_n}}{\rd m_{d_n}^{(n)}} = c_{d_n,\kappa_n}^{(n)} \int_{\mathrm{SO}(d_n+1)} \prod_{i=1}^n \exp(\kappa_n\, \bX_{n,i}^\top \bO \btheta) \,\rd\bO,\label{eq:ivMF}
\end{align}
where $c_{d_n,\kappa_n}^{(n)}$ is a normalizing constant and the integral is with respect to the Haar measure on~$\mathrm{SO}(d_n+1)$. Unlike the von Mises--Fisher distribution with density $\bx \mapsto c_{d_n,\kappa_n} \exp(\kappa_n\, \bx^\top\btheta)$, the integrated von Mises--Fisher distribution is invariant with respect to rotations. This invariance justifies its choice for analyzing the power of uniformity tests, as these are obviously invariant with respect to rotations and thus their asymptotic power cannot depend on $\btheta$. As shown in \citet[Theorem 4.2]{Cutting2017}, the Rayleigh test enjoys non-trivial asymptotic power against such integrated von Mises--Fisher alternatives with $\kappa_n \asymp d_n^{3/4}/\sqrt{n}$.

The following result gives the high-dimensional asymptotic power of Sobolev tests against alternatives where \eqref{eq:ivMF} holds, extending that in \citet[Section 4.2]{Cutting2017} for the Rayleigh~test.

\begin{theorem} \label{thm:LANstuff}
Assume that the setup of Theorem \ref{thm:asymp} and its assumptions \ref{cond1} and \ref{cond2} hold for the situation in which $\bX_{n,1},\ldots,\bX_{n,n}$ are distributed as integrated von Mises--Fisher on $\Sdn$ with concentration $\kappa_n = \tau_n d_n^{3/4}/\sqrt{n}$ for the positive sequence~$\tau_n\to\tau>0$. Assume moreover the limit
\begin{align}
    \Gamma\defin \lim_{n \to \infty} \frac{\sqrt{d_n} {v}^2_{1,d_n}}{\sqrt{2}\big(\sum_{k=1}^\infty v_{k,d_n}^4 d_{k,d_n}\big)^{1/2}}\label{eq:Gamma}
\end{align}
exists.

Then $\sigma_n^{-1}T_n\inlaw\mathcal{N}(\Gamma \tau^2,1)$.
\end{theorem}

As a consequence of Theorem \ref{thm:LANstuff}, the one-sided test that rejects at the asymptotic level $\alpha$ when $\sigma_n^{-1}T_n>z_\alpha$, where $z_\alpha$ is the upper $\alpha$-quantile of a standard normal, has asymptotic power $1-\Phi(z_{\alpha}-\Gamma\tau^2)$ against contiguous alternatives given by integrated von Mises--Fisher distributions with $\kappa_n = \tau_n d_n^{3/4}/\sqrt{n}$. Therefore, it is able to detect such alternatives if and only if $\Gamma>0$. Since $d_{k_0,d_n}\sim d_n^{k_0}/k_0!$, we have two relevant implications:
\begin{itemize}
    \item For the subclass of \emph{$k_0$-Sobolev tests} with $v_{k,d_n}=\delta_{k,k_0}$, $\Gamma=(1/\sqrt{2}) 1_{\{k_0=1\}}$. Therefore, the Rayleigh test is the only one of these tests detecting the studied contiguous alternatives. In particular, the Bingham test with $v_{k,d_n}=\delta_{k,2}$ is ``blind'' against them.
    \item For the subclass of \emph{finite Sobolev tests} with $v_{k,d_n}=1_{\{k\leq k_0\}}$, for $k_0>1$, $\Gamma=0$. Therefore, none of these tests can detect the studied contiguous alternatives, despite the similarity between the Rayleigh sequence $v_{k,d_n}=\delta_{k,1}$ and, e.g., the ``hybrid Rayleigh--Bingham'' test arising with $v_{k,d_n}=\delta_{k,1}+\delta_{k,2}$.
\end{itemize}

It is possible to combine the Rayleigh test with other Sobolev tests so that the resulting test is still able to detect contiguous integrated von Mises--Fisher alternatives. This can be achieved by making the sequence of weights $(v_{k,d_n})_{k=1}^\infty$ decay with $d_n$. Indeed, it is immediate to check that, given $v_{1,d_n}=1$, weights of the form ${v}_{k,d_n}\sim \big[a_{k,d_n} k!\, d_n^{-(k-1)}\big]^{1/4}$ for $k>1$ with $A\defin\lim_{n\to\infty}\sum_{k=2}^{\infty}a_{k,d_n}$ existing in $[0,\infty)$ yield $\Gamma=1/\sqrt{2(1+A)}$.

Appendix C
in the SM reports numerical experiments illustrating Theorems \ref{thm:asymp} and \ref{thm:LANstuff} and the effect of the above-discussed sequences of weights on the local power.

\section{Applications to testing problems}
\label{sec:ext}

\subsection{Rotational symmetry}
\label{sec:rotasym}

A random vector~$\bX$ on~$\Sdn$ has a \emph{rotationally symmetric} distribution about $\btheta \in \Sdn$ if and only if~$\bO \bX$ has the same distribution as $\bX$ for any orthogonal matrix~$\bO\in\mathrm{SO}(d_n+1)$ such that $\bO\btheta=\btheta$. Rotational symmetry is closely related to the \emph{tangent-normal decomposition}: for any~$\bx\in \Sdn$,
\begin{align}
    \bx = v_{\btheta}(\bx) \btheta +
    \big(1-v^2_{\btheta}(\bx)\big)^{1/2} \, {\bGamma}_{\btheta} \bu_{\btheta}(\bx), \quad
    v_{\btheta}(\bx)
    	\defin
    	\bx^\top\btheta,\quad
    	\bu_{\btheta}(\bx)
    	\defin
    	\frac{{\bGamma}_{\btheta}^\top \bx}{\|{\bGamma}_{\btheta}^\top\bx\|}, \label{eq:tangnorm}
\end{align}
where ${\bGamma}_{\btheta}$ denotes an arbitrary $(d_n+1) \times d_n$ matrix whose columns form an orthogonal complement to~$\btheta$ (so that~${\bGamma}_{\btheta}^\top{\bGamma}_{\btheta}=\bI_{d_n}$ and ${\bGamma}_{\btheta} {\bGamma}_{\btheta}^\top=\bI_{d_n+1} - \btheta \btheta^\top$). If $\bX$ is rotationally symmetric about $\btheta$, then the \emph{multivariate sign}~$\bu_{\btheta}(\bX)$ is uniformly distributed over~$\Sdnm$. Therefore, based on a sample $\bX_{n,1}, \ldots, \bX_{n,n}$ on $\Sdn$, one way to test for rotational symmetry about $\btheta$ is to test the necessary condition of the uniformity of $\bu_{\btheta}(\bX_{n,1}), \ldots, \bu_{\btheta}(\bX_{n,n})$ on $\Sdnm$.

The tangent von Mises--Fisher distribution defined in \citet[Section 2.2]{Garcia-Portugues2020} provides a set of alternatives to rotational symmetry about $\btheta$ based on \eqref{eq:tangnorm}. A random vector $\bX$ on $\Sdn$ has such a distribution if $v_{\btheta}(\bX)$ follows an arbitrary density $g$ on $[-1,1]$, $\bu_{\btheta}(\bX)$ is von Mises--Fisher distributed with location $\bmu\in \Sdnm$ and concentration $\kappa_n\geq 0$ (with $\kappa_n=0$ corresponding to rotational symmetry about $\btheta$), and $v_{\btheta}(\bX)$ and $\bu_{\btheta}(\bX)$ are independent. As in \eqref{eq:ivMF}, an \emph{integrated tangent von Mises--Fisher distribution} is characterized by the likelihood function
\begin{align}
    \frac{\rd{\rm P}^{(n)}_{\btheta,g,\kappa_n}}{\rd m_{d_n}^{(n)}} = c_{d_n,g,\kappa_n}^{(n)}  \prod_{i=1}^n g(v_{\btheta}(\bX_{n,i}))\int_{\mathrm{SO}(d_n)} \prod_{i=1}^n \exp(\kappa_n\, \bu_{\btheta}(\bX_{n,i})^\top \bO \bmu) \,\rd\bO.\label{eq:liktvMF}
\end{align}

We have the following result as a direct consequence of Theorem \ref{thm:LANstuff}, \eqref{eq:ivMF}, and \eqref{eq:liktvMF}.

\begin{corollary}
Let $d_n$ be a sequence of positive integers such that $d_n\to\infty$ as $n\to\infty$. Assume that $\bX_{n,i}$, $i=1,\ldots,n$, $n=1,2,\ldots$, form a triangular array such that, for any fixed $n$, the random vectors $\bX_{n,1},\ldots, \bX_{n,n}$ are distributed as integrated tangent von Mises--Fisher about $\btheta\in\Sdn$ with concentration $\kappa_n = \tau_n d_n^{3/4}/\sqrt{n}$ for the positive sequence~$\tau_n\to\tau>0$. Assume moreover that assumptions \ref{cond1} and \ref{cond2} of Theorem \ref{thm:asymp} hold for the sample $\bu_{\btheta}(\bX_{n,1}), \ldots, \bu_{\btheta}(\bX_{n,n})$, and that $T_n$ is computed on this sample. Then $\sigma_n^{-1}T_n\inlaw\mathcal{N}(\Gamma\tau^2,1)$ with $\Gamma$ given in \eqref{eq:Gamma}.
\end{corollary}

\subsection{Spherical symmetry}
\label{sec:sphsym}

A random vector~$\bX$ on~$\R^{d_n+1}$ has a \emph{spherically symmetric} distribution if and only if~$\bO \bX$ has the same distribution as $\bX$ for any orthogonal matrix~$\bO\in\mathrm{SO}(d_n+1)$. In the absolutely continuous case, spherical symmetry of $\bX$ is characterized by the projection $\pi(\bX)\defin\bX/\|\bX\|$ being uniformly distributed over $\Sdn$ and the radius $\|\bX\|$ being independent of $\pi(\bX)$. Therefore, based on a sample $\bX_{n,1}, \ldots, \bX_{n,n}$ on $\R^{d_n+1}$, one way to test for spherical symmetry is to test the necessary condition of the uniformity of $\pi(\bX_{n,1}), \ldots, \pi(\bX_{n,n})$ on $\Sdn$, where Theorem \ref{thm:asymp} is applicable. As explained in \cite{Garcia-Portugues:Sobolev}, Sobolev tests can therefore be seen as multivariate sign tests of spherical symmetry. Spherical symmetry has many applications, including genetics \citep{Chen2010tests}, finance \citep{yang2021testing}, and astrophysics (Section \ref{sec:cmb}).

The previous characterization of spherical symmetry can be exploited to construct high-dimensional goodness-of-fit tests for distributions on $\R^{d_n+1}$, i.e., test
\begin{align}
    \Hcal_{0,n}: \bX_{n,1}, \ldots, \bX_{n,n}\text{ share the spherically symmetric distribution } F_{0,n} \label{eq:gof}
\end{align}
for a specified distribution $F_{0,n}$ on $\R^{d_n+1}$. For that purpose, we only need a classical (one-dimensional) goodness-of-fit test for the common radial distribution $R_{0,n}$ of the sample $\|\bX_{n,1}\|,\ldots, \|\bX_{n,n}\|$ that is induced by $F_{0,n}$. Let $A_n$ denote a test statistic for this purpose, such as the Anderson--Darling test statistic, and let $F_A(x)\defin\lim_{n\to\infty}\Pr[A_n\leq x]$, $x\geq0$, be the null asymptotic distribution of $A_n$. Then, given the sample $\bX_{n,1}, \ldots, \bX_{n,n}$ on $\R^{d_n+1}$, we test \eqref{eq:gof} by combining tests of (\textit{i}) that $\|\bX_{n,1}\|, \ldots, \|\bX_{n,n}\|$ share the radial distribution $R_{0,n}$ and (\textit{ii}) that $\pi(\bX_{n,1}), \ldots, \pi(\bX_{n,n})$ are uniformly distributed on $\Sdn$. Under \eqref{eq:gof}, the radial and projected samples are independent, so the corresponding $p$-values are asymptotically independent and can be aggregated by the Fisher method \citep{Fisher1925}
\begin{align}
    G_n\defin-2\big[\log(1-F_A(A_n))+\log\big(1-F_1\big((\sigma_n^{-1}T_n)^2\big)\big)\big], \label{eq:twosided}
\end{align}
where $F_1$ denotes the distribution function of a $\chi^2_1$, $A_n$ is computed with $\|\bX_{n,1}\|, \ldots, \|\bX_{n,n}\|$, and $T_n$ is computed with $\pi(\bX_{n,1}), \ldots, \pi(\bX_{n,n})$. Under the null hypothesis in \eqref{eq:gof}, as $d_n\to\infty$ with $n\to\infty$, the asymptotic distribution of $G_n$ is $\chi^2_4$. Recall that \eqref{eq:twosided} features a two-sided version of the test of uniformity on $\Sdn$; an alternative one-sided version is $G_n\defin -2 [\log(1-F_A(A_n))+\log(1-\Phi(\sigma_n^{-1}T_n))]$. Some particularly relevant cases for the distribution $F_{0,n}$ in \eqref{eq:gof} are:
\begin{itemize}
    \item Multivariate standard normal $\mathcal{N}_{d_n+1}(\zero,\bI_{d_n+1})$, where $\|\bX\|^2$ follows a $\chi^2_{d_n+1}$ distribution.
    \item Multivariate standard Student's $t$ distribution with $\nu>0$ degrees of freedom, where $\|\bX\|^2/(d_n+1)$ follows Snedecor's $F_{d_n+1,\nu}$ distribution.
    \item Multivariate isotropic stable distribution, where $\|\bX\|^2$ is distributed as the product $AT$, where $A$ follows a positive stable law $\mathrm{S}(\beta/2,1,2\gamma_0^2(\cos(\pi\beta/4))^{2/\beta},0)$, $0<\beta<2$, $\gamma_0>0$, and $T$ independently follows a $\chi^2_{d_n+1}$ distribution \citep{Nolan2013}.
\end{itemize}

The simple null hypothesis tested in \eqref{eq:gof} can be generalized to the composite null
\begin{align}
    \Hcal_{0,n}: \bX_{n,1}, \ldots, \bX_{n,n}\text{ share the spherically symmetric distribution } F_{0,n,\btheta} \label{eq:goftheta}
\end{align}
with $\btheta\in\Theta\subset\R^q$ an unknown parameter only affecting the radial distribution $R_{0,n,\btheta}$ of $\|\bX_{n,1}\|, \ldots,\allowbreak \|\bX_{n,n}\|$. In this case, the above testing pipeline only has to be modified in the radial part to carry out a parametric goodness-of-fit test, with $A_n$ suitably adapted and $F_A$ potentially approximated by parametric bootstrapping. A relevant case of \eqref{eq:goftheta} is the high-dimensional goodness-of-fit for spherically symmetric distributions on $\R^{d_n+1}$ with squared radius distributed as $\Gamma(\theta_1,\theta_2)$, which in particular comprises standard normality for $\theta_1=(d_n+1)/2$ and $\theta_2=2$. Another interesting case is the multivariate $t$ distribution where the degrees of freedom $\nu$ are unknown and have to be estimated.

\section{Simulation study}
\label{sec:num}

We present the results of a quadripartite Monte Carlo simulation study conducted to demonstrate the power performance of the goodness-of-fit test for spherical symmetry presented in Section \ref{sec:sphsym} in high-dimensional settings. For all four parts we used the following simulation setting. To test the radial component, we utilized the classical Anderson--Darling test for \eqref{eq:gof} and incorporated parameter estimation for \eqref{eq:goftheta}. In the latter (composite) case we employed a bootstrap algorithm which is implemented as follows:
\begin{enumerate}
    \item Compute the estimator $\hat{\boldsymbol{\theta}}$ from the radial part of the data.
    \item Compute the Anderson--Darling test statistic $A_n$ assuming $\hat{\boldsymbol{\theta}}$ is the underlying parameter.
    \item Generate $B$ bootstrap samples from $R_{0,n,\hat{\boldsymbol{\theta}}}$, the estimated radial distribution (recall \eqref{eq:goftheta}).
    \item Compute for each bootstrap sample $j$ the estimator $\tilde{\boldsymbol{\theta}}_j$ and compute $\tilde{A}_{n,j}$ assuming $\tilde{\boldsymbol{\theta}}_j$ is the underlying parameter, $j=1,\ldots,B$.
    \item Generate the bootstrap $p$-value $p_{A_n}=\frac{1}{B}\sum_{j=1}^B1_{\{\tilde{A}_{n,j}> A_n\}}$.
\end{enumerate}
In the final computation of $G_n$, the directional part of the data is kept fixed, and the bootstrap $p$-value $p_{A_n}$ replaces the $p$-value of the radial part.

We conducted $M=5{,}000$ Monte Carlo repetitions with a fixed significance level of $5\pct$. When employing a bootstrap algorithm, we used $B=500$ bootstrap samples. To construct alternatives to \eqref{eq:gof} and \eqref{eq:goftheta}, we used the following models: $t_\nu(\bS)$, the centered multivariate $t$ distribution with $\nu > 0$ degrees of freedom and scale matrix $\bS$; $st_\nu(\bS,\bxi)$, the skewed version of $t_\nu(\bS)$ by the skewing vector $\bxi$ \citep{A:2014}; $\mathrm{vMF}(\kappa)$, the von Mises--Fisher distribution with mean direction $(1, 0, \ldots, 0)^\top$ and concentration $\kappa > 0$; $F \otimes G$, the law of $\bX \cdot Y$, where $\bX$ is a random vector with distribution $F$ and $Y$ is a positive random variable with distribution $G$. We set $\bxi_1 = (1, 0, \ldots, 0)^\top $ and take $\bOmega_1$ to be the matrix with ones on the diagonal and $0.25$ in every off-diagonal entry. For the multivariate normal distribution, the covariance matrix is denoted by $\bSigma_p =\mathrm{diag}(0.5, \stackrel{\lfloor pd \rfloor}{\ldots}, 0.5,\allowbreak 1.5, \stackrel{d - \lfloor pd \rfloor}{\ldots}, 1.5)$. The identity matrix of size $d\times d$ is denoted by $\bI_d$. The model $\mathrm{DMN}(\rho)$ (dependent multivariate normal), with $\rho \in [0, 1/\sqrt{2}]$, is defined by $\bX\cdot R$, where: $R=[F^{-1}_{d}(\Phi(Z))]^{1/2}$, with $\Phi$ and $F_d$ being the distribution functions of a standard normal and a $\chi^2_{d}$, respectively; $\bX=\bY/\|\bY\|$; and $(Z,\bY^\top)^\top$ is distributed as $\mathcal{N}_{d+1}(\zero,\bSigma_\rho)$ for $\bSigma_\rho=(1, \brho^\top; \brho, \bI_{d})$ and $\brho=(\rho^1,\ldots,\rho^d)^\top$. Note that $R$ and $\bX$ have the marginal distributions of the radial and projected components of $\mathcal{N}_{d}(\zero,\bI_{d})$, but that dependence between $R$ and $\bX$ is present if $\rho>0$. We denote by $|F|$ the distribution of $|X|$ when $X$ is distributed as $F$, $\mathrm{Ca}(\mu, \gamma)$ refers to the Cauchy distribution with location parameter $\mu\in\mathbb{R}$ and scale parameter $\gamma>0$, and $\Gamma(k,\theta)$ stands for a gamma distribution with shape $k>0$ and scale $\theta>0$. Lastly, $\mathrm{S}(\beta) = \mathrm{S}(\beta/2, 1, 2\gamma_0^2(\cos(\pi\beta/4))^{2/\beta}, 0)$, for $0<\beta < 2$ and $\gamma_0=1$, denotes a positive stable law.

In the first part of the simulation study, we provide results designed to demonstrate the power performance of the goodness-of-fit test of high-dimensional normality. As a benchmark, we selected the simple hypothesis version of the classical BHEP test for multivariate normality, which is based on the weighted $L^2$ distance between the characteristic function of the standard normal distribution and its empirical counterpart, as described in \cite{EH:2020}. This test is equivalent to the one used in \cite{HW:1997}, except that it is applied to the original data instead of scaled residuals. Critical values for the BHEP test, depending on the sample size $n$ and dimension $d$, were obtained via Monte Carlo simulation with $100{,}000$ repetitions. The results are presented in Table \ref{tab:norm}. We observe that both tests are well-calibrated, although the BHEP test appears conservative when $d > n$. In lower dimensions (e.g., $d = 100$), the BHEP test outperforms the new test for most alternatives, except for the multivariate Student's $t$ distribution. However, as the dimensionality increases ($d = 200, 300$), the new test significantly surpasses the BHEP test. Notably, for the $t$ distribution, as the degrees of freedom increase and approach the null hypothesis, the new test retains power, while the BHEP test completely fails.

\begin{table}[t]
\centering
\caption{\small Empirical rejection rates (in percentage) for testing high-dimensional normality in comparison to the simple hypothesis version of the BHEP test of multivariate normality.}\label{tab:norm}
\setlength{\tabcolsep}{5.25pt}
\begin{tabular}{r|rrr|rrr||rrr|rrr}
\toprule
  $n$ & \multicolumn{6}{c||}{$100$} & \multicolumn{6}{c}{$200$} \\
$d$ & $100$ & $200$ & $300$ & $100$ & $200$ & $300$ & $100$ & $200$ & $300$ & $100$ & $200$ & $300$ \\
  Test & \multicolumn{3}{c|}{$G_n$} & \multicolumn{3}{c||}{BHEP} & \multicolumn{3}{c|}{$G_n$} & \multicolumn{3}{c}{BHEP} \\
  \midrule
$\mathcal{N}_{d}(\zero,\bI_d)$ & 4 & 5 & 5 & 5 & 1 & 0 & 5 & 5 & 5 & 6 & 2 & 1 \\
$\mathcal{N}_{d}(\zero,\bSigma_{0.1})$ & 6 & 5 & 6 & 8 & 2 & 1 & 6 & 6 & 6 & 9 & 2 & 1 \\
$\mathcal{N}_{d}(\zero,\bSigma_{0.25})$ & 8 & 8 & 8 & 12 & 3 & 1 & 9 & 9 & 8 & 20 & 6 & 2 \\
$\mathcal{N}_{d}(\zero,\bSigma_{0.5})$ & 13 & 14 & 14 & 22 & 7 & 2 & 19 & 20 & 20 & 49 & 17 & 6 \\
$\mathrm{vMF}(1)\otimes [\chi^2_{d}]^{1/2}$ & 6 & 5 & 5 & 6 & 1 & 0 & 6 & 5 & 5 & 7 & 2 & 0 \\
$\mathrm{vMF}(4)\otimes [\chi^2_{d}]^{1/2}$ & 18 & 7 & 5 & 29 & 4 & 1 & 47 & 11 & 7 & 65 & 7 & 1 \\
$t_{10}(0.75\cdot \bI_d)$ & 100 & 100 & 100 & 49 & 48 & 48 & 100 & 100 & 100 & 86 & 88 & 88 \\
$t_{10}(0.8\cdot \bI_d)$ & 100 & 100 & 100 & 14 & 11 & 10 & 100 & 100 & 100 & 31 & 30 & 29 \\
$t_{10}(0.9\cdot \bI_d)$ & 100 & 100 & 100 & 20 & 17 & 17 & 100 & 100 & 100 & 28 & 25 & 25 \\
$t_{10}(1.1\cdot \bI_d)$ & 100 & 100 & 100 & 100 & 100 & 100 & 100 & 100 & 100 & 100 & 100 & 100 \\
$t_{10}(\bI_d)$ & 100 & 100 & 100 & 84 & 87 & 87 & 100 & 100 & 100 & 98 & 99 & 99 \\
$t_{30}(\bI_d)$ & 100 & 100 & 100 & 23 & 17 & 17 & 100 & 100 & 100 & 34 & 30 & 31 \\
$t_{100}(\bI_d)$ & 86 & 100 & 100 & 7 & 2 & 1 & 99 & 100 & 100 & 8 & 3 & 1 \\
$t_{500}(\bI_d)$ & 9 & 21 & 43 & 5 & 1 & 0 & 13 & 42 & 76 & 6 & 2 & 0 \\
$t_{1000}(\bI_d)$ & 6 & 9 & 14 & 5 & 2 & 0 & 7 & 12 & 24 & 5 & 1 & 0 \\
\bottomrule
\end{tabular}
\end{table}

The second part of the simulation study is two-fold. First, we test the simple null hypothesis $\mathcal{H}_{0,n}$ as defined in \eqref{eq:gof}, which posits that the data follow a multivariate standard Student's $t$ distribution with fixed degrees of freedom $\nu$. Second, we examine the composite hypothesis in \eqref{eq:goftheta}, wherein the parameter $\nu$ is estimated using the maximum-likelihood estimator with respect to the radial component of the data. In the finite-dimensional case, \cite{MMOV:2024} investigated the goodness-of-fit testing of the family of multivariate $t$ distributions. Notably, all known methods provided or cited therein exhibit severe computational limitations as the dimension $d$ increases, in some cases as early as $d\geq5$. Therefore, to the best of our knowledge, the methodology presented herein is unrivaled.

\begin{table}[t]
\centering
\caption{\small Empirical rejection rates (in percentage) for testing the fit to the high-dimensional multivariate Student distribution with known $\nu=5$ degrees of freedom (top) and unknown degrees of freedom parameter estimated by maximum likelihood (bottom).}\label{tab:Stud}
\setlength{\tabcolsep}{4.6pt}
\begin{tabular}{r|rrr|rrr|rrr|rrr}
\toprule
  $n$ & \multicolumn{3}{c}{$100$} & \multicolumn{3}{c}{$200$} & \multicolumn{3}{c}{$500$} & \multicolumn{3}{c}{$1000$} \\
 $d$ & $100$ & $200$ & $300$ & $100$ & $200$ & $300$ & $100$ & $200$ & $300$ & $100$ & $200$ & $300$ \\
  \midrule
$t_5(\bI_d)$ & 5 & 5 & 5 & 5 & 5 & 5 & 5 & 5 & 4 & 5 & 5 & 5 \\
$t_6(\bI_d)$ & 6 & 6 & 6 & 9 & 10 & 10 & 23 & 24 & 22 & 53 & 55 & 56 \\
$t_7(\bI_d)$ & 12 & 13 & 13 & 27 & 29 & 31 & 81 & 84 & 85 & 100 & 100 & 100 \\
$t_8(\bI_d)$ & 24 & 24 & 25 & 61 & 64 & 65 & 100 & 100 & 100 & 100 & 100 & 100 \\
$t_5(0.9\cdot \bI_d)$ & 26 & 26 & 27 & 46 & 47 & 47 & 86 & 87 & 89 & 99 & 100 & 100 \\
$\mathrm{vMF}(5)\otimes[d\cdot$F$_{d,5}]^{1/2}$ & 34 & 9 & 6 & 81 & 20 & 10 & 100 & 75 & 33 & 100 & 100 & 80 \\
$\mathrm{vMF}(10)\otimes[d\cdot$F$_{d,5}]^{1/2}$ & 100 & 57 & 22 & 100 & 97 & 62 & 100 & 100 & 100 & 100 & 100 & 100 \\
$st_5(\bOmega_1,\zero)$ & 29 & 40 & 49 & 38 & 50 & 57 & 65 & 73 & 78 & 88 & 93 & 94 \\
$st_5(\bI_d,\bxi_1)$ & 48 & 28 & 20 & 95 & 76 & 60 & 100 & 100 & 100 & 100 & 100 & 100 \\
$st_8(\bOmega_1,\zero)$ & 43 & 55 & 61 & 68 & 78 & 81 & 98 & 99 & 99 & 100 & 100 & 100 \\
$st_8(\bI_d,\bxi_1)$ & 71 & 54 & 45 & 100 & 97 & 95 & 100 & 100 & 100 & 100 & 100 & 100 \\
\midrule
$t_5(\bI_d)$ & 5 & 5 & 5 & 5 & 5 & 5 & 6 & 6 & 5 & 5 & 5 & 5 \\
$t_5(0.8\cdot \bI_d)$ & 86 & 87 & 87 & 99 & 99 & 99 & 100 & 100 & 100 & 100 & 100 & 100 \\
$t_5(0.9\cdot \bI_d)$ & 28 & 30 & 30 & 51 & 51 & 52 & 88 & 89 & 90 & 100 & 100 & 100 \\
$t_5(1.1\cdot \bI_d)$ & 19 & 20 & 20 & 38 & 39 & 40 & 79 & 80 & 80 & 98 & 98 & 99 \\
$t_5(1.25\cdot \bI_d)$ & 79 & 81 & 81 & 99 & 99 & 99 & 100 & 100 & 100 & 100 & 100 & 100 \\
$\mathrm{vMF}(1)\otimes[d\cdot$F$_{d,5}]^{1/2}$ & 5 & 5 & 5 & 6 & 5 & 5 & 7 & 5 & 6 & 11 & 6 & 6 \\
$\mathrm{vMF}(3)\otimes[d\cdot$F$_{d,5}]^{1/2}$ & 9 & 5 & 6 & 23 & 8 & 7 & 73 & 17 & 9 & 100 & 48 & 21 \\
$\mathrm{vMF}(5)\otimes[d\cdot$F$_{d,5}]^{1/2}$ & 33 & 8 & 7 & 80 & 20 & 10 & 100 & 73 & 31 & 100 & 100 & 81 \\
$\mathrm{vMF}(10)\otimes[d\cdot$F$_{d,5}]^{1/2}$ & 100 & 56 & 24 & 100 & 97 & 64 & 100 & 100 & 100 & 100 & 100 & 100 \\
$st_5(\bOmega_1,\zero)$ & 29 & 41 & 50 & 38 & 49 & 57 & 63 & 71 & 77 & 86 & 91 & 93 \\
$st_5(\bI_d,\bxi_1)$ & 48 & 29 & 21 & 95 & 76 & 61 & 100 & 100 & 100 & 100 & 100 & 100 \\
\bottomrule
\end{tabular}
\end{table}

The results in Table \ref{tab:Stud} present the empirical power performance of both hypotheses. In both cases we see that the method is well-calibrated. For the simple hypothesis case, in the upper half of the table, we detect changes in the parameter $\nu$ as well as disturbances in the scale matrix. Notably, under these alternatives the power is stable across different dimensions, since the disturbance only affects the radial part of the distribution. Disturbances in the directional part are also well identified. Dependence and asymmetry are modeled in the multivariate skewed $t$ alternative and are both well detected, although the power of the tests decreases in the asymmetric case when the dimension increases. For the composite case, in the lower half of the table, the procedure effectively identifies substantial disturbances in the scale matrix; however, small anomalies in the directional component prove to be more challenging to detect. The power behavior under the multivariate skewed $t$ distribution interestingly parallels the simple hypothesis case.

The third part of the simulation study consists of testing for a high-dimensional multivariate stable law for which, to the best of our knowledge, no competing methods exist in the literature. In Table \ref{tab:stable}, we present the simulation results for testing the multivariate stable distribution for $\beta=1$ in dimensions $d=50$ and $d=100$. Accordingly, the null hypothesis corresponds to $\mathcal{N}_{d}(\zero,\bI_d)\otimes[\mathrm{S}(1)]^{1/2}$. It is important to note that there is no numerically stable algorithm for this class of distributions in higher dimensions; see \cite{Nolan2013}. The test is able to detect deviations in the context of the multivariate $t$ distribution as the degrees of freedom increase but shows little sensitivity to alternatives with dependence structures, such as DMN. However, when the radial part is also changed to another stable distribution, the power of the test increases, notably in higher dimensions.

\begin{table}[!htb]
\centering
\caption{\small Empirical rejection rates (in percentage) for testing the high-dimensional multivariate stable distribution with $\beta=1$, here $\mathcal{N}_{d}(\zero,\bI_d)\otimes[\mathrm{S}(1)]^{1/2}$.}\label{tab:stable}
\begin{tabular}{r|rr|rr}
\toprule
$n$ & \multicolumn{2}{c|}{$50$} & \multicolumn{2}{c}{$100$} \\
$d$ & $50$ & $100$ & $50$ & $100$ \\
  \midrule
$\mathcal{N}_{d}(\zero,\bI_d)\otimes[\mathrm{S}(1)]^{1/2}$ & 5 & 5 & 5 & 5 \\
$t_1(\bI_d)$ & 5 & 5 & 5 & 5 \\
$t_{1.25}(\bI_d)$ & 8 & 8 & 13 & 13 \\
$t_{1.5}(\bI_d)$ & 17 & 17 & 41 & 42 \\
$\mathrm{DMN}(0)\otimes[\mathrm{S}(0.8)]^{1/2}$ & 39 & 43 & 69 & 75 \\
$\mathrm{DMN}(0.25)\otimes[\mathrm{S}(0.8)]^{1/2}$ & 40 & 43 & 69 & 75 \\
$\mathrm{DMN}(0.5)\otimes[\mathrm{S}(0.8)]^{1/2}$ & 40 & 44 & 70 & 75 \\
\bottomrule
\end{tabular}
\end{table}

Fourth and finally, we present results for the composite case of testing for high-dimensional spherically symmetric distributions with radii following an unknown gamma distribution. The results in Table \ref{tab:gamma.est} demonstrate that the procedure, using the bootstrap algorithm with maximum-likelihood estimators for both parameters of the gamma distribution, is well-calibrated. It effectively detects disturbances in the directional component as the concentration parameter of the vMF distribution increases. When examining changes in the radial distribution with only minor disturbances in the directional component, the power remains consistent across all dimensions, as expected, since the power of the test predominantly arises from the radial component.

\begin{table}[!htb]
\centering
\caption{\small Empirical rejection rates (in percentage) for testing the fit to the high-dimensional spherically symmetric distribution with radial part from the gamma family of distributions, both parameters estimated by maximum likelihood.}\label{tab:gamma.est}
\begin{tabular}{r|rrrr|rrrr}
\toprule
$n$ & \multicolumn{4}{c|}{$100$} & \multicolumn{4}{c}{$200$} \\
$d$ & $100$ & $200$ & $300$ & $1000$ & $100$ & $200$ & $300$ & $1000$ \\
  \midrule
$\mathrm{vMF}(0)\otimes\Gamma(2,5)$ & 5 & 4 & 5 & 6 & 5 & 5 & 6 & 5 \\
$\mathrm{vMF}(0.25)\otimes\chi^2_2$ & 5 & 5 & 4 & 6 & 5 & 5 & 5 & 6 \\
$\mathrm{vMF}(5)\otimes\chi^2_{d}$ & 33 & 10 & 6 & 5 & 80 & 22 & 9 & 5 \\
$\mathrm{vMF}(10)\otimes\chi^2_{20}$ & 100 & 56 & 21 & 5 & 100 & 98 & 61 & 8 \\
$\mathrm{vMF}(2)\otimes\Gamma(2,5)$ & 6 & 6 & 5 & 5 & 9 & 6 & 6 & 6 \\
$\mathrm{vMF}(5)\otimes\Gamma(2,5)$ & 34 & 10 & 6 & 5 & 80 & 20 & 9 & 6 \\
$\mathrm{vMF}(10)\otimes\Gamma(2,5)$ & 100 & 55 & 23 & 5 & 100 & 98 & 62 & 7 \\
$\mathrm{vMF}(20)\otimes\Gamma(2,5)$ & 100 & 100 & 99 & 11 & 100 & 100 & 100 & 32 \\
$\mathrm{vMF}(0.25)\otimes|\mathrm{Ca}(2,5)|$ & 98 & 98 & 98 & 98 & 100 & 100 & 100 & 100 \\
$\mathrm{vMF}(0.5)\otimes|t_2|$ & 48 & 48 & 48 & 48 & 75 & 75 & 75 & 75 \\
\bottomrule
\end{tabular}
\end{table}

\section{Gaussianity of the cosmic microwave background}
\label{sec:cmb}

The Cosmic Microwave Background (CMB) is the relic radiation released about $380{,}000$ years after the Big Bang. Under the standard cosmological model, once its mean level and the motion-induced dipole are removed, the CMB temperature on the celestial sphere is a single realization of a statistically isotropic Gaussian random field \citep{Planck2015_XVI}; assessing this prediction is a basic test of the model. We consider the foreground-cleaned SMICA temperature map of the Planck mission \citep{Planck2015_IX}, read within R through the \texttt{rcosmo} package \citep{rcosmo} on the HEALPix equal-area pixelization of the sphere \citep{Gorski2005}, to investigate the isotropic-Gaussian field hypothesis.

A random field $Z(\bt)$ on $\mathbb{S}^2$ is summarized by its spherical-harmonic coefficients $a_{\ell m}$:
$$
Z(\bt) = \sum_{\ell=0}^\infty\sum_{m=-\ell}^\ell a_{\ell m}Y_{\ell m}(\bt),\qquad \bt\in\mathbb{S}^2,
$$
where $Y_{\ell m}$ are the degree-$\ell$ spherical harmonics on $\mathbb{S}^2$ (unlike \eqref{Harmo}, we adopt here the standard notation in physics for the indices). If $Z(\bt)$ is isotropic and Gaussian, the $a_{\ell m}$ are independent zero-mean normal variables whose variance $C_\ell$ depends only on $\ell$. Dividing the $2\ell+1$ coefficients of each degree by their empirical root-mean-square therefore yields a pool of approximately independent standard normal variables, and arranging $n(d_n+1)$ of them into an $n\times(d_n+1)$ matrix produces vectors $\bX_{n,1},\ldots,\bX_{n,n}$ on $\R^{d_n+1}$ that, under the hypothesis, are approximately independent $\mathcal{N}_{d_n+1}(\zero,\bI_{d_n+1})$, with projections $\pi(\bX_{n,i})$ uniform on $\Sdn$. Testing the isotropic-Gaussian hypothesis thus amounts to testing the multivariate standard normality of $\bX_{n,1},\ldots,\bX_{n,n}$.

We coarsen the full-sky map to the $N_{\mathrm{side}}=64$ HEALPix grid (Figure \ref{fig:cmb}) and compute the coefficients $a_{\ell m}$ for $\ell=2,\ldots,\ell_{\max}$ by quadrature over the grid, excluding the monopole $\ell=0$ and dipole $\ell=1$ as is standard in CMB analyses. We consider $\ell_{\max}=100$ as it stays well within the resolution $\ell<2N_{\mathrm{side}}$ of the grid. We set $n=d_n=100$ for the $n\times(d_n+1)$ matrix, as this choice (\emph{i}) lies in the high-dimensional regime $n,d_n\to\infty$ of our asymptotic theory; and (\emph{ii}) accommodates $n(d_n+1)=10{,}100$ coefficients, hence maximizing the use of the available coefficients $\sum_{\ell=2}^{100}(2\ell+1)=10{,}197$.

\begin{figure}[!htb]
\centering
\includegraphics[width=0.85\textwidth,trim={0cm 1.1cm 0cm 1.1cm},clip]{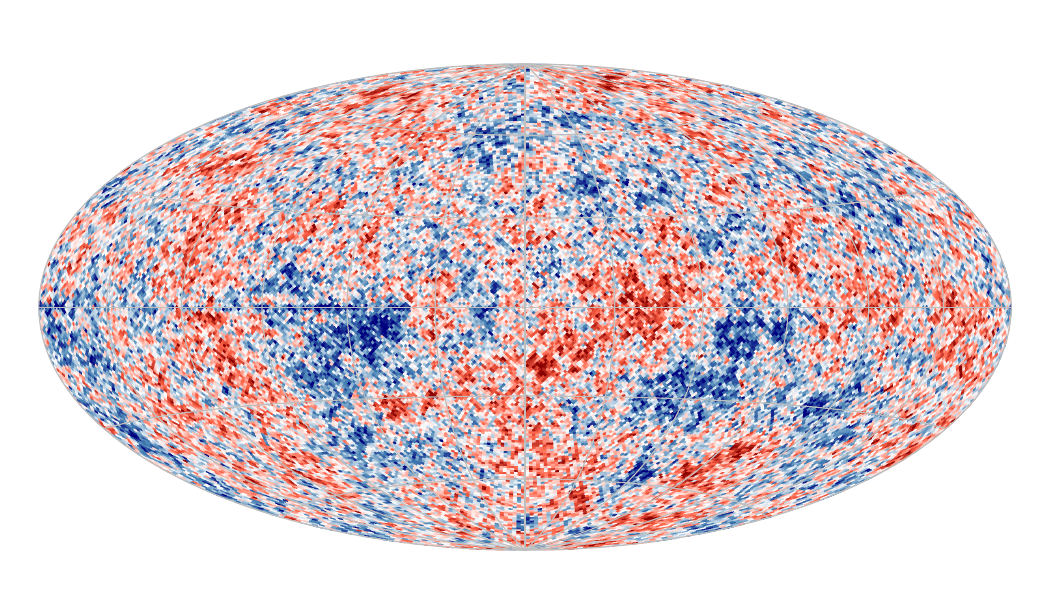}
\caption{\small Planck SMICA CMB temperature map (Hammer projection) on the $N_{\mathrm{side}}=64$ HEALPix grid used in the analysis. The first $10{,}100$ standardized spherical-harmonic coefficients of this map, reshaped into $n=100$ vectors on $\R^{101}$, are consistent with high-dimensional uniformity and Gaussianity (Table \ref{tab:cmb}).}\label{fig:cmb}
\end{figure}

We test the uniformity of these directions with the $k_0$-Sobolev test statistics in \eqref{eq:asympCkexact} for $k_0=1,\ldots,6$, using the asymptotic null distribution and considering both a one-sided (upper-tail, $\mathcal{N}(0,1)$ null) and a two-sided ($\chi^2_1$ null) version of the test. The one-sided version is natural against concentration alternatives, while the two-sided version also flags departures towards over-regularity. For each $k_0$ we also perform the associated test of high-dimensional normality, aggregating by Fisher's method (Section \ref{sec:sphsym}) the $k_0$-Sobolev test with an Anderson--Darling test of the fit of $\|\bX_{n,i}\|^2$ to $\chi^2_{101}$; the resulting aggregated statistic is asymptotically $\chi^2_4$ (recall \eqref{eq:twosided}).

\begin{table}[!htb]
\centering
\caption{\small Asymptotic $p$-values for standardized spherical harmonic coefficients of the Planck SMICA CMB map. The coefficients are reshaped into $n=100$ vectors on $\R^{101}$. The first two rows give the $k_0$-Sobolev tests of uniformity of the projections on $\mathbb{S}^{100}$, in one- and two-sided form; the last two rows give the associated tests of high-dimensional normality, aggregating with Fisher's method an Anderson--Darling test of the radial component.}\label{tab:cmb}
\begin{tabular}{l|c|cccccc}
\toprule
Test & Sided & $k_0=1$ & $k_0=2$ & $k_0=3$ & $k_0=4$ & $k_0=5$ & $k_0=6$ \\
\midrule
Uniformity & One & $0.679$ & $0.961$ & $0.129$ & $0.706$ & $0.854$ & $0.844$ \\
 & Two & $0.642$ & $0.078$ & $0.259$ & $0.589$ & $0.292$ & $0.311$ \\
Normality & One & $0.929$ & $0.996$ & $0.381$ & $0.939$ & $0.981$ & $0.979$ \\
 & Two & $0.912$ & $0.267$ & $0.591$ & $0.885$ & $0.633$ & $0.657$ \\
\bottomrule
\end{tabular}
\end{table}

Table \ref{tab:cmb} reports the results of the tests. None of the tests rejects at conventional significance levels. The one-sided Rayleigh test yields a $p$-value of $0.68$, and the one-sided $k_0$-Sobolev $p$-values range from $0.13$ to $0.96$; their two-sided counterparts range from $0.08$ to $0.64$. The shared radial Anderson--Darling test of $\|\bX_{n,i}\|^2\sim\chi^2_{101}$ gives $p$-value $0.95$, and the Fisher combinations remain far from rejection. The directions extracted from the Planck map are thus statistically indistinguishable from uniform, and the map is consistent with the isotropic Gaussian prediction of the standard cosmological model.

\section{Discussion}
\label{sec:disc}

We have derived the asymptotic null distribution for the entire class of Sobolev tests of uniformity on the hypersphere, accommodating arbitrary growth rates of both the dimension and sample size. Additionally, we have established the consistency and characterized the limiting behavior of these tests under sequences of integrated von Mises--Fisher alternatives. Our findings have been successfully applied to high-dimensional testing problems involving rotational and spherical symmetry. These results demonstrate that the proposed testing procedures are well-calibrated, maintaining statistical power against both local and fixed alternatives.

The present contribution opens a line of future research on the exhaustive determination of the high-dimensional detection thresholds of Sobolev tests of uniformity against general sequences of rotationally symmetric alternatives, paralleling the developments in the fixed-dimensional case in \cite{Garcia-Portugues:Sobolev}.

\section*{Supplementary materials}

The supplementary materials contain numerical experiments illustrating Theorems \ref{thm:asymp}--\ref{thm:LANstuff}.

\section*{Acknowledgments}

The second author acknowledges support from grant PID2021-124051NB-I00, funded by MCIN/\-AEI/\-10.13039/\-501100011033 and by ``ERDF A way of making Europe''. His research was also supported by ``Convocatoria de la Universidad Carlos III de Madrid de Ayudas para la recualificación del sistema universitario español para 2021--2023'', funded by Spain's Ministerio de Ciencia, Innovación y Universidades. The third author acknowledges support from PDR and Hubert-Curien grants from the Fonds National de la Recherche Scientifique (FNRS). Comments by two anonymous reviewers are acknowledged.

\appendix

\section{Main proofs}
\label{app:proofs}

\begin{proof}[Proof of Theorem \ref{thm:asymp}]
As in \citet[Theorem 2.1]{Paindaveine2016}, the proof is based on an application of \cite{Billingsley2012}'s (central limit) Theorem 35.12 for martingale difference sequences. For that reason, the theorem is stated as Theorem \ref{thm:billingsley} in Appendix \ref{app:lemmas}.

The proof that $\sigma_n^{-1} T_n\inlaw\mathcal{N}(0,1)$ is split into checking the three main conditions of Theorem \ref{thm:billingsley}: that $T_n$ has a martingale difference structure, that the variances of the martingale difference terms are finite, and that the Lindeberg condition is satisfied.

\smallskip
\noindent\textit{Martingale difference structure.}
\smallskip

We consider the $\sigma$-algebra $\mathcal{F}_{n,\ell}=\sigma(\{\bX_{n,1},\ldots,\bX_{n,\ell}\})$ that is spanned by the random vectors $\bX_{n,1},\ldots,\bX_{n,\ell}$ on $\Sdn$.

It is beneficial to define standardized versions of the kernel $\psi_n$ and the statistic $T_n$ as
$$
\psint(x)\defin\sigma_n^{-1}\psi_n(x)=:\sum_{k=1}^\infty \lrp{1+\frac{2k}{d_n-1}}\tilde{v}^2_{k,d_n} C_k^{(d_n-1)/2}(x)\quad\text{and}\quad\tilde{T}_n\defin\frac{2}{n}\sum_{1\leq i<j\leq n} \psint(\bX_{n,i}^\top\bX_{n,j}).
$$
Under the null hypothesis of uniformity (assumed henceforth), the kernel $\psint$ verifies that
\begin{align}
    \E{\psint(\bX_{n,1}^\top\bX_{n,2})}&=0\quad\text{and} \label{eq:zeroexp1}\\
    \E{\psint(\bX_{n,1}^\top\bX_{n,2})\mid \bX_{n,1}}&=0. \label{eq:zeroexp2}
\end{align}
Indeed, \eqref{eq:zeroexp1} directly follows from \eqref{eq:zeroexp2} and the latter from
\begin{align*}
    \E{\psint(\bX_{n,1}^\top\bX_{n,2})\mid \bX_{n,1}}&=\int_{\Sdn} \psint(\bX_{n,1}^\top\bx_2)\,\rd\nu_{d_n}(\bx_2)\\
    &=\frac{1}{\om{d_n}}\int_{\Sdnm}\int_{-1}^1 \psint(x)(1-x^2)^{d_n/2-1}\,\rd x\,\rd\nu_{d_n-1}(\bxi)\\
    &=\frac{\om{d_n-1}}{\om{d_n}}\int_{-1}^1 \sum_{k=1}^\infty \lrp{1+\frac{2k}{d_n-1}}\tilde{v}^2_{k,d_n}C_k^{(d_n-1)/2}(x)(1-x^2)^{d_n/2-1}\,\rd x\\
    &=0,
\end{align*}
where $\nu_{d_n}$ denotes the uniform measure on $\Sdn$. Above, we have applied, first, the tangent-normal change of variables $\bx=x\bmu+(1-x^2)^{1/2}\bB_{\bmu} \bxi$, for $\bx,\bmu\in\Sdn$, $x=\bx^\top\bmu$, $\bxi\in\Sdnm$, and $\bB_{\bmu}$ a $(d_n+1)\times d_n$ matrix such that $\bB_{\bmu}^\top\bB_{\bmu}=\bI_{d_n}$ and $\bB_{\bmu}\bB_{\bmu}^\top=\bI_{d_n+1}-\bmu\bmu^\top$; and then the orthogonality of the Gegenbauer polynomials ($C_0^{(d_n-1)/2}\equiv1$).

We define
\begin{align}
    D_{n,\ell}\defin\E{\tilde{T}_n\mid \mathcal{F}_{n,\ell}}-\E{\tilde{T}_n\mid \mathcal{F}_{n,\ell-1}},\quad \ell=1,\ldots,n,\, n=1,2,\ldots,\label{eq:Dnl}
\end{align}
which, by definition, is a martingale difference sequence. The form of $D_{n,\ell}$ is determined by computing $\E{\tilde{T}_n\mid \mathcal{F}_{n,\ell}}$ next. Throughout, sums over empty sets of indices are defined as being equal to zero. We have that
\begin{align*}
    \frac{n}{2}\E{\tilde{T}_n\mid \mathcal{F}_{n,\ell}}=&\;\sum_{i=1}^{n-1}\sum_{j=i+1}^n \E{\psint(\bX_{n,i}^\top\bX_{n,j})\mid \mathcal{F}_{n,\ell}}\\
    =&\; \sum_{i=1}^{\ell-1}\sum_{j=i+1}^{\ell-1} \E{\psint(\bX_{n,i}^\top\bX_{n,j})\mid \bX_{n,i},\bX_{n,j}}\\
    &+\sum_{i=1}^{\ell-1} \E{\psint(\bX_{n,i}^\top\bX_{n,\ell})\mid \bX_{n,i},\bX_{n,\ell}}+\sum_{i=1}^{\ell-1}\sum_{j=\ell+1}^{n-1} \E{\psint(\bX_{n,i}^\top\bX_{n,j})\mid \bX_{n,i}}\\
    &+ \sum_{j=\ell+1}^n \E{\psint(\bX_{n,\ell}^\top\bX_{n,j})\mid \bX_{n,\ell}}+ \sum_{i=\ell+1}^{n-1}\sum_{j=i+1}^n \E{\psint(\bX_{n,i}^\top\bX_{n,j})\mid \bX_{n,1},\ldots,\bX_{n,\ell}}\\
    &=\sum_{i=1}^{\ell-1}\sum_{j=i+1}^{\ell-1} \psint(\bX_{n,i}^\top\bX_{n,j})+\sum_{i=1}^{\ell-1} \psint(\bX_{n,i}^\top\bX_{n,\ell})\\
    &=\sum_{i=1}^{\ell-1}\sum_{j=i+1}^{\ell} \psint(\bX_{n,i}^\top\bX_{n,j}),
\end{align*}
where we have used \eqref{eq:zeroexp1} and \eqref{eq:zeroexp2}. Therefore,
\begin{align*}
\frac{n}{2}D_{n,\ell}&=\frac{n}{2}\lrp{\E{\tilde{T}_n\mid \mathcal{F}_{n,\ell}}-\E{\tilde{T}_n\mid \mathcal{F}_{n,\ell-1}}}\\
    &=\sum_{i=1}^{\ell-1}\sum_{j=i+1}^\ell \psint(\bX_{n,i}^\top\bX_{n,j})-\sum_{i=1}^{\ell-2}\sum_{j=i+1}^{\ell-1} \psint(\bX_{n,i}^\top\bX_{n,j})\\
    &=\psint(\bX_{n,\ell-1}^\top\bX_{n,\ell})+\sum_{i=1}^{\ell-2}\sum_{j=i+1}^\ell \psint(\bX_{n,i}^\top\bX_{n,j})-\sum_{i=1}^{\ell-2}\sum_{j=i+1}^{\ell-1} \psint(\bX_{n,i}^\top\bX_{n,j})\\
    &=\sum_{i=1}^{\ell-1} \psint(\bX_{n,i}^\top\bX_{n,\ell}).
\end{align*}
In particular, we have that $D_{n,1}=0$. Hence,
\begin{align*}
    \tilde{T}_n&=\frac{2}{n}\sum_{1\leq i<j\leq n} \psint(\bX_{n,i}^\top\bX_{n,j})=\frac{2}{n}\sum_{j=2}^n\sum_{i=1}^{j-1} \psint(\bX_{n,i}^\top\bX_{n,j})=\sum_{j=2}^n D_{n,j}=\sum_{\ell=1}^n D_{n,\ell}
\end{align*}
and $\tilde{T}_n$ has a martingale difference structure.

\smallskip
\noindent\textit{Finite variance of $D_{n,\ell}$.}
\smallskip

First, observe that $\psint(\bX_{n,i}^\top\bX_{n,\ell})$ and $\psint(\bX_{n,j}^\top\bX_{n,\ell})$ are uncorrelated for $i\neq j$ and $1\leq i,j<\ell$ since
\begin{align*}
    \Ebig{\psint(\bX_{n,i}^\top\bX_{n,\ell}) \psint(\bX_{n,j}^\top\bX_{n,\ell})}=&\;\EBig{\Ebig{\psint(\bX_{n,i}^\top\bX_{n,\ell})\mid\bX_{n,\ell}} \Ebig{\psint(\bX_{n,j}^\top\bX_{n,\ell})\mid\bX_{n,\ell}}}=0
\end{align*}
because of \eqref{eq:zeroexp2}. We then have that
\begin{align*}
    \V{D_{n,\ell}}=&\;\frac{4}{n^2}\sum_{i=1}^{\ell-1}\V{\psint(\bX_{n,i}^\top\bX_{n,\ell})}
    =\frac{4}{n^2}(\ell-1)\E{\psint^2(\bX_{n,1}^\top\bX_{n,2})}.
\end{align*}
Because of the definition of $\sigma_n^2$ in \eqref{eq:sigman2},
\begin{align}
    \E{\psint^2(\bX_{n,1}^\top\bX_{n,2})}=&\;\frac{\Ebig{\psin^2(\bX_{n,1}^\top\bX_{n,2})}}{\sigma_n^2}=\frac{1}{2} \label{eq:Epsi2}
\end{align}
and hence
\begin{align}
    \V{D_{n,\ell}}=&\;\frac{2}{n^2}(\ell-1)<\infty.\label{eq:varDnl}
\end{align}

\smallskip
\noindent\textit{Convergence of $\sum_{\ell=1}^{n} \sigma^2_{n,\ell}$ to $1$ in probability.}
\smallskip

We show that $\sum_{\ell=1}^{n} \sigma^2_{n,\ell}$ converges to $1$ in mean square, where $\sigma^2_{n,\ell}\defin\Ebig{D_{n,\ell}^{2} \mid \mathcal{F}_{n,\ell-1}}$. For that, we show that $\Ebig{\sum_{\ell=1}^{n} \sigma^2_{n,\ell}}\to 1$ and $\Vbig{\sum_{\ell=1}^{n} \sigma^2_{n,\ell}}\to 0$.

First, note that
\begin{align}
    \sigma^2_{n,\ell}=&\;\E{D_{n,\ell}^{2} \mid \mathcal{F}_{n,\ell-1}}\nonumber\\
    =&\;\frac{4}{n^2}\sum_{i=1}^{\ell-1} \E{\psint^2(\bX_{n,i}^\top\bX_{n,\ell})\mid \bX_{n,i}}+\frac{4}{n^2}\sum_{1\leq i\neq j\leq \ell-1} \Ebig{\psint(\bX_{n,i}^\top\bX_{n,\ell})\psint(\bX_{n,j}^\top\bX_{n,\ell})\mid \bX_{n,i},\bX_{n,j}}\nonumber\\
    =&\;\frac{4}{n^2} \frac{\ell-1}{2}+\frac{8}{n^2}\sum_{1\leq i< j\leq \ell-1} \Ebig{\psint(\bX_{n,i}^\top\bX_{n,\ell})\psint(\bX_{n,j}^\top\bX_{n,\ell})\mid \bX_{n,i},\bX_{n,j}}\label{eq:sigma2nl}
\end{align}
because
\begin{align*}
    \E{\psint^2(\bX_{n,i}^\top\bX_{n,\ell})\mid \bX_{n,i}}=\frac{\om{d_n-1}}{\om{d_n}}\int_{-1}^1 \psint^2(x)(1-x^2)^{d_n/2-1}\,\rd x=\E{\psint^2(\bX_{n,i}^\top\bX_{n,\ell})}=\frac{1}{2}.
\end{align*}

Then, from \eqref{eq:sigma2nl} and \eqref{eq:Epsi2},
\begin{align*}
    \E{\sigma^2_{n,\ell}}=\frac{2(\ell-1)}{n^2}
\end{align*}
since, for $i\neq j$ and $1\leq i,j\leq \ell-1$, $\Ebig{\psint(\bX_{n,i}^\top\bX_{n,\ell})\psint(\bX_{n,j}^\top\bX_{n,\ell})}=\Ebig{\Ebig{\psint(\bX_{n,i}^\top\bX_{n,\ell})\mid \bX_{n,\ell}}\cdot \allowbreak \Ebig{\psint(\bX_{n,j}^\top\bX_{n,\ell})\mid \bX_{n,\ell}}}=0$ by \eqref{eq:zeroexp2}. It therefore follows that
\begin{align*}
    \E{\sum_{\ell=1}^n\sigma^2_{n,\ell}}=\frac{2}{n^2}\sum_{\ell=1}^n(\ell-1)=\frac{2}{n^2}\sum_{\ell=1}^{n-1}\ell=\frac{2}{n^2}\frac{(n-1)n}{2}=\frac{n-1}{n}\to 1
\end{align*}
as $n\to\infty$.

Now, from \eqref{eq:sigma2nl}, we have that
\begin{align}
    \V{\sum_{\ell=1}^n\sigma^2_{n,\ell}}
    &=\frac{64}{n^4}\V{\sum_{\ell=3}^n\sum_{1\leq i<j\leq \ell-1} \E{\psint(\bX_{n,i}^\top\bX_{n,\ell})\psint(\bX_{n,j}^\top\bX_{n,\ell})\mid \bX_{n,i},\bX_{n,j}}}\nonumber\\
    &=\frac{64}{n^4}\V{\sum_{1\leq i<j\leq n} (n-j)\E{\psint(\bX_{n,i}^\top\bX_{n,\ell})\psint(\bX_{n,j}^\top\bX_{n,\ell})\mid \bX_{n,i},\bX_{n,j}}}\label{eq:cijsum}\\
    &\indef\frac{64}{n^4}\V{\sum_{1\leq i<j\leq n} (n-j)f_{ij}},\label{eq:fij}
\end{align}
where in passing to \eqref{eq:cijsum} we have used that $\sum_{\ell=3}^n \sum_{1\leq i< j\leq \ell-1} c_{ij}=\sum_{1\leq i<j\leq n-1} (n-j)c_{ij}=\sum_{1\leq i<j\leq n} (n-j) c_{ij}$. In \eqref{eq:cijsum} and what follows, $\ell$ represents any index different from $i$ and $j$. The $\{f_{ij}\}_{1\leq i<j\leq n}$ are pairwise uncorrelated, since
\begin{align*}
    \E{f_{12}\mid\bX_{n,1}}&=\E{\psint(\bX_{n,1}^\top\bX_{n,\ell})\psint(\bX_{n,2}^\top\bX_{n,\ell})\mid \bX_{n,1}}\\
    &=\E{\psint(\bX_{n,1}^\top\bX_{n,\ell})\E{\psint(\bX_{n,2}^\top\bX_{n,\ell})\mid \bX_{n,\ell}}\mid \bX_{n,1}}=0
\end{align*}
and therefore $\E{f_{12}f_{13}}=\E{\E{f_{12}\mid\bX_{n,1}}\E{f_{13}\mid\bX_{n,1}}}=0$ and $\E{f_{12}}=0$. Hence, \eqref{eq:fij} becomes
\begin{align*}
    \V{\sum_{\ell=1}^n\sigma^2_{n,\ell}}=&\;\frac{64}{n^4}\sum_{1\leq i<j\leq n} (n-j)^2\V{f_{ij}}\\
    =&\;\frac{64}{n^4}\E{f_{12}^2}\sum_{1\leq i<j\leq n} (n-j)^2.
\end{align*}
Now,
\begin{align*}
\sum_{1\leq i<j\leq n} (n-j)^2=\sum_{j=2}^n (j-1)(n-j)^2=\sum_{j=1}^{n-1} j(n-j-1)^2\leq n^2\sum_{j=1}^{n-1} j=\frac{n^3(n-1)}{2},
\end{align*}
which gives
\begin{align*}
    \V{\sum_{\ell=1}^n\sigma^2_{n,\ell}}\leq\frac{64}{n^4}\E{f_{12}^2}\frac{n^3(n-1)}{2}=\frac{32(n-1)}{n}\E{f_{12}^2}.
\end{align*}

Exploiting the definition \eqref{eq:psi} of the kernel $\psint$ in terms of a mean square convergent series, we can express $\Ebig{f_{ij}^2}$ through the coefficients of $\psin$:
\begin{align}
    \mathrm{E}\Big[\psint&(\bX_{n,i}^\top\bX_{n,\ell})\psint(\bX_{n,j}^\top\bX_{n,\ell})\mid \bX_{n,i},\bX_{n,j}\Big]\nonumber\\
    =&\; \int_{\Sdn} \psint(\bX_{n,i}^\top\bx)\psint(\bX_{n,j}^\top\bx)\,\rd\nu_{d_n}(\bx) \nonumber\\
    =&\sum_{k_1,k_2=1}^\infty \lrp{1+\frac{2k_1}{d_n-1}}\lrp{1+\frac{2k_2}{d_n-1}}\tilde{v}^2_{k_1,d_n}\tilde{v}^2_{k_2,d_n}\int_{\Sdn}\!\! C_{k_1}^{(d_n-1)/2}(\bX_{n,i}^\top\bx)C_{k_2}^{(d_n-1)/2}(\bX_{n,j}^\top\bx)\,\rd\nu_{d_n}(\bx) \nonumber\\
    =&\; \sum_{k=1}^\infty \lrp{1+\frac{2k}{d_n-1}}^{2} \tilde{v}^4_{k,d_n} \lrp{1+\frac{2k}{d_n-1}}^{-1}C_{k}^{(d_n-1)/2}(\bX_{n,i}^\top\bX_{n,j})\label{eq:CkXiXj}\\
    =&\; \sum_{k=1}^\infty \lrp{1+\frac{2k}{d_n-1}} \tilde{v}^4_{k,d_n} C_{k}^{(d_n-1)/2}(\bX_{n,i}^\top\bX_{n,j}),\nonumber
\end{align}
with equality \eqref{eq:CkXiXj} following, e.g., from Lemma B.7 in the Supplementary Material of \cite{Garcia-Portugues2023}. Then, due to the orthogonality of Gegenbauer polynomials,
\begin{align}
    \Ebig{f_{ij}^2}=&\;\E{\lrb{ \sum_{k=1}^\infty \lrp{1+\frac{2k}{d_n-1}} \tilde{v}^4_{k,d_n} C_{k}^{(d_n-1)/2}(\bX_{n,i}^\top\bX_{n,j}) }^2}\nonumber\\
    =&\;\E{\sum_{k=1}^\infty\lrp{1+\frac{2k}{d_n-1}}^2 \tilde{v}^8_{k,d_n} \lrp{C_{k}^{(d_n-1)/2}(\bX_{n,i}^\top\bX_{n,j})}^2}\nonumber\\
    =&\;\sum_{k=1}^\infty\lrp{1+\frac{2k}{d_n-1}}^2 \tilde{v}^8_{k,d_n} \frac{\om{d_n-1}}{\om{d_n}}c_{k,d_n}, \label{eq:Efij2}
\end{align}
with $c_{k,d_n}$ defined in \eqref{eq:ckd1}. Due to \eqref{eq:ckd2}, \eqref{eq:Efij2} becomes
\begin{align*}
    \Ebig{f_{ij}^2}=\sum_{k=1}^\infty \tilde{v}^8_{k,d_n} d_{k,d_n}=\frac{\sum_{k=1}^\infty v^8_{k,d_n} d_{k,d_n}}{4\Ebig{\psin^2(\bX_{n,1}^\top\bX_{n,2})}^2}=\frac{\sum_{k=1}^\infty v^8_{k,d_n} d_{k,d_n}}{4\lrpbig{\sum_{k=1}^\infty v^4_{k,d_n} d_{k,d_n}}^2}.
\end{align*}
Therefore,
\begin{align*}
    \V{\sum_{\ell=1}^n\sigma^2_{n,\ell}}\leq\frac{32(n-1)}{n}\E{f_{12}^2}=\frac{8(n-1)}{n} \frac{\sum_{k=1}^\infty v^8_{k,d_n} d_{k,d_n}}{\lrpbig{\sum_{k=1}^\infty v^4_{k,d_n} d_{k,d_n}}^2}=o(1)
\end{align*}
by condition \eqref{eq:Av8}.

\smallskip
\noindent\textit{Lindeberg condition.}
\smallskip

By the Cauchy--Schwarz and Chebyshev inequalities, and then \eqref{eq:varDnl}, we have
\begin{align}
    \sum_{\ell=1}^{n} \E{D_{n,\ell}^2 1_{\{|D_{n,\ell}|>\varepsilon\}}}\leq&\; \sum_{\ell=1}^{n} \sqrt{\Ebig{D_{n,\ell}^4}} \sqrt{\Prob{|D_{n,\ell}|>\varepsilon}} \nonumber\\
    \leq &\; \frac{1}{\varepsilon}\sum_{\ell=1}^{n} \sqrt{\Ebig{D_{n,\ell}^4}} \sqrt{\V{D_{n,\ell}}} \nonumber\\
    =&\; \frac{\sqrt{2}}{\varepsilon n}\sum_{\ell=1}^{n} \sqrt{\Ebig{D_{n,\ell}^4}} \sqrt{\ell-1}. \label{eq:lind1}
\end{align}

The fourth moment of $D_{n,\ell}$ is
\begin{align}
    \E{D_{n,\ell}^4}=&\;\frac{16}{n^4}\E{\lrp{\sum_{i=1}^{\ell-1} \psint(\bX_{n,i}^\top\bX_{n,\ell})}^4}\nonumber\\
    =&\;\frac{16}{n^4}\sum_{i_1,i_2,i_3,i_4=1}^{\ell-1}\E{\psint(\bX_{n,i_1}^\top\bX_{n,\ell})\psint(\bX_{n,i_2}^\top\bX_{n,\ell})\psint(\bX_{n,i_3}^\top\bX_{n,\ell})\psint(\bX_{n,i_4}^\top\bX_{n,\ell})}\nonumber\\
    =&\;\frac{16}{n^4}\lrb{(\ell-1)\E{\psint^4(\bX_{n,1}^\top\bX_{n,\ell})}+3(\ell-1)(\ell-2)\E{\psint^2(\bX_{n,1}^\top\bX_{n,\ell})}^2}\label{eq:twoterms}\\
    =&\;\frac{16(\ell-1)}{n^4}\lrb{\E{\psint^4(\bX_{n,1}^\top\bX_{n,\ell})}+\frac{3}{4}(\ell-2)}\nonumber
\end{align}
with \eqref{eq:twoterms} following because of \eqref{eq:zeroexp1} and \eqref{eq:zeroexp2}, which makes other terms equal to zero (e.g., due to \eqref{eq:zeroexp2}, $\Ebig{\psint(\bX_{n,1}^\top\bX_{n,\ell})\psint^3(\bX_{n,2}^\top\bX_{n,\ell})}=0$). Since $\Ebig{\psint^4(\bX_{n,1}^\top\bX_{n,\ell})}=O(1)$ as $n\to\infty$ and $d_n\to\infty$ by condition \eqref{eq:Av8},
\begin{align}
    \E{D_{n,\ell}^4}\leq &\;\frac{16(\ell-1)^2}{n^4}\lrb{C+3}\label{eq:Dnl4}
\end{align}
for $C>0$ independent of $n$.

Substituting \eqref{eq:Dnl4} into \eqref{eq:lind1} results in
\begin{align*}
    \sum_{\ell=1}^{n} \E{D_{n,\ell}^2 1_{\{|D_{n,\ell}|>\varepsilon\}}} \leq &\; \frac{\sqrt{2}}{\varepsilon n}\sum_{\ell=1}^{n} \sqrt{\frac{16(\ell-1)^2(C+3)}{n^4}} \sqrt{\ell-1}\\
    = &\; \frac{4\sqrt{2(C+3)}}{\varepsilon n^3}\sum_{\ell=1}^{n} (\ell-1)^{3/2}\\
    = &\; \frac{4\sqrt{2(C+3)}}{\varepsilon n^3}O(n^{5/2})\\
    = &\; \frac{1}{\varepsilon}O(n^{-1/2}),
\end{align*}
and hence the Lindeberg condition is also satisfied.
\end{proof}

\begin{proof}[Proof of Theorem \ref{thm:consis}]
The proof of the theorem is an almost direct consequence of the following lemma.
\begin{lemma} \label{thm:asympaltgen} Assume that $\bX_{n,1},\ldots, \bX_{n,n}$ are mutually independent and identically distributed with common distribution $F_n$. Then,
$\mathrm{E}_{F_n}[T_{n}]=(n-1) \sum_{k=1}^\infty v_{k,d_n}^2 \sum_{r=1}^{d_{k,d_n}} e_{r,k,n}^2$ and
$$
\mathrm{Var}_{F_n}[T_{n}]= \frac{2 (n-1)}{n} \sum_{k, k'=1}^\infty v_{k,d_n}^2 v_{k',d_n}^2 \sum_{r=1}^{d_{k,d_n}} \sum_{r'=1}^{d_{k',d_n}} \big\{2(n-1) e_{r,k,n} e_{r',k',n} c_{r,r',n}^{(k,k')}+ \big(c_{r,r',n}^{(k,k')}\big)^2 \big\}.
$$
\end{lemma}
\begin{proof}[Proof of Lemma \ref{thm:asympaltgen}]
It follows directly from Lemma A in \citet[p. 183]{Serfling1980} that
\begin{align}
    \V{T_n}
    &= \mathrm{Var}_{F_n}\Big[(n-1) \frac{2}{n(n-1)} \sum_{1 \leq i < j \leq n} \psin(\bX_{n,i}^\top\bX_{n,j})\Big]\nonumber\\
    &= \frac{2 (n-1)}{n} (2 (n-2) \xi_1+ \xi_2), \label{eq:varserf}
\end{align}
where
\begin{align} \label{xi1}
    \xi_1&=\mathrm{E}_{F_n}\big[(\mathrm{E}_{F_n}[\psin(\bX_{n,1}^\top\bX_{n,2})\mid\bX_{n,2}])^2\big]-\mathrm{E}_{F_n}^2[\psin(\bX_{n,1}^\top\bX_{n,2})]
\end{align}
and
\begin{align} \label{xi2}
    \xi_2&=\mathrm{Var}_{F_n}[\psin(\bX_{n,1}^\top\bX_{n,2})].
\end{align}
Using \eqref{Harmo}, we easily obtain that
\begin{align*}
    \mathrm{E}_{F_n}\big[\psin(\bX_{n,1}^\top\bX_{n,2})\big]&=\sum_{k=1}^\infty v_{k,d_n}^2 \sum_{r=1}^{d_{k,d_n}} e_{r,k,n}^2,\\
    \mathrm{E}_{F_n}\big[\psin(\bX_{n,1}^\top\bX_{n,2})\mid\bX_{n,2}\big]&= \sum_{k=1}^\infty v_{k,d_n}^2 \sum_{r=1}^{d_{k,d_n}} e_{r,k,n} \; g_{r, k}(\bX_{n,2})
\end{align*}
and that
\begin{align*}
    \mathrm{E}_{F_n}\big[\psin^2(\bX_{n,1}^\top\bX_{n,2})\big] &= \sum_{k,k'=1}^\infty v_{k,d_n}^2 v_{k',d_n}^2 \sum_{r=1}^{d_{k,d_n}} \sum_{r'=1}^{d_{k',d_n}}  \mathrm{E}_{F_n}^2[g_{r, k}(\bX_{n,1})g_{r', k'}(\bX_{n,1})].
\end{align*}

Plugging these three expressions into \eqref{xi1} and \eqref{xi2} gives
\begin{align*}
    \xi_1
    =&\;\sum_{k,k'=1}^\infty v_{k,d_n}^2 v_{k',d_n}^2 \sum_{r=1}^{d_{k,d_n}} \sum_{r'=1}^{d_{k',d_n}}\Big\{e_{r,k,n}\, e_{r',k',n}\, \mathrm{E}_{F_n}[g_{r,k}(\bX_{n,2})\, g_{r',k'}(\bX_{n,2})]\\
    &- e_{r,k,n}^2\, e_{r',k',n}^2\Big\}\\
    =&\;\sum_{k,k'=1}^\infty v_{k,d_n}^2 v_{k',d_n}^2 \sum_{r=1}^{d_{k,d_n}} \sum_{r'=1}^{d_{k',d_n}} e_{r,k,n}\, e_{r',k',n}\, c_{r,r',n}^{(k,k')}
\end{align*}
and
\begin{align*}
    \xi_2
    =&\;\sum_{k,k'=1}^\infty v_{k,d_n}^2 v_{k',d_n}^2 \sum_{r=1}^{d_{k,d_n}} \sum_{r'=1}^{d_{k',d_n}} \Big\{\mathrm{E}_{F_n}^2[g_{r,k}(\bX_{n,1})\, g_{r',k'}(\bX_{n,1})] - e_{r,k,n}^2\, e_{r',k',n}^2\Big\}\\
    =&\;\sum_{k,k'=1}^\infty v_{k,d_n}^2 v_{k',d_n}^2 \sum_{r=1}^{d_{k,d_n}} \sum_{r'=1}^{d_{k',d_n}} \Big\{\big(c_{r,r',n}^{(k,k')}\big)^2 + 2\, c_{r,r',n}^{(k,k')}\, e_{r,k,n}\, e_{r',k',n}\Big\}.
\end{align*}

Replacing into \eqref{eq:varserf} yields the stated result:
\begin{align*}
\mathrm{Var}_{F_n}[T_n]
    =&\;\frac{2(n-1)}{n}\Bigg\{2(n-2)\sum_{k,k'=1}^\infty v_{k,d_n}^2 v_{k',d_n}^2 \sum_{r=1}^{d_{k,d_n}} \sum_{r'=1}^{d_{k',d_n}} e_{r,k,n}\, e_{r',k',n}\, c_{r,r',n}^{(k,k')} \\
    &+ \sum_{k,k'=1}^\infty v_{k,d_n}^2 v_{k',d_n}^2 \sum_{r=1}^{d_{k,d_n}} \sum_{r'=1}^{d_{k',d_n}} \Big[\big(c_{r,r',n}^{(k,k')}\big)^2 + 2\, c_{r,r',n}^{(k,k')}\, e_{r,k,n}\, e_{r',k',n}\Big]\Bigg\} \\
    =&\;\frac{2(n-1)}{n}\sum_{k,k'=1}^\infty v_{k,d_n}^2 v_{k',d_n}^2 \sum_{r=1}^{d_{k,d_n}} \sum_{r'=1}^{d_{k',d_n}} \Big\{2(n-1)\, e_{r,k,n}\, e_{r',k',n}\, c_{r,r',n}^{(k,k')} + \big(c_{r,r',n}^{(k,k')}\big)^2\Big\}.
\end{align*}
\end{proof}

It directly follows from Lemma \ref{thm:asympaltgen} above that
$$
\lim_{n \to \infty} \mathrm{E}_{F_n}[\sigma_n^{-1} n^{-1} T_{n}]=\lim_{n \to \infty}  \sigma_n^{-1} \sum_{k=1}^\infty v_{k,d_n}^2 \sum_{r=1}^{d_{k,d_n}} e_{r,k,n}^2=c
$$
thanks to condition \ref{thm:consis:1} and that
\begin{align*}
    \lim_{n \to \infty} \mathrm{Var}_{F_n}[\sigma_n^{-1} n^{-1} T_{n}]
    =&\;\lim_{n \to \infty} 4 \sigma_n^{-2} n^{-1} \sum_{k, k'=1}^\infty v_{k,d_n}^2 v_{k',d_n}^2 \sum_{r=1}^{d_{k,d_n}} \sum_{r'=1}^{d_{k',d_n}}  e_{r,k,n} e_{r',k',n} c_{r,r',n}^{(k,k')} \\
    & + \lim_{n \to \infty} 2\sigma_n^{-2} n^{-2} \sum_{k, k'=1}^\infty v_{k,d_n}^2 v_{k',d_n}^2 \sum_{r=1}^{d_{k,d_n}} \sum_{r'=1}^{d_{k',d_n}}\big(c_{r,r',n}^{(k,k')}\big)^2 \\
=&\;o(1)
\end{align*}
as $n \to \infty$ thanks to conditions \ref{thm:consis:2} and \ref{thm:consis:3}. As a direct consequence, $\sigma_n^{-1} n^{-1} T_{n}$ converges in probability to the strictly positive constant $c$ and therefore $\sigma_n^{-1}  T_{n}$ diverges to $\infty$ in probability as $n \to \infty$ and
$$\lim_{n \to \infty} \mathrm{E}_{F_n}[\phi_{n}]=1,$$
which is the desired result.
\end{proof}

\begin{proof}[Proof of Theorem \ref{thm:LANstuff}] It directly follows from Theorem 4.2 in \cite{Cutting2017} that the sequence of integrated von Mises--Fisher models with concentration $\kappa_n$ is locally asymptotically normal with central sequence
$$
\sum_{\ell=1}^n \Delta_{n,\ell}\defin \sum_{\ell=1}^n \frac{\sqrt{d_n}}{n} \sum_{i=1}^{\ell-1} \bX_{n,i}^\top \bX_{n,\ell}
$$
and Fisher information $\Sigma= 1/2$.

Under the null hypothesis, both $D_{n,\ell}$ (defined in \eqref{eq:Dnl} in the proof of Theorem \ref{thm:asymp}) and $\Delta_{n,\ell}$ are martingale differences with respect to the filtration $\mathcal{F}_{n,\ell}=\sigma(\{\bX_{n,1},\ldots,\bX_{n,\ell}\})$. Hence, for every $a,b\in\R$, $M_{n,\ell}(a,b)\defin aD_{n,\ell}+b\Delta_{n,\ell}$ is a martingale difference array satisfying the assumptions of Theorem \ref{thm:billingsley}: $\sum_{\ell=1}^n\E{M_{n,\ell}(a,b)^2\mid\mathcal{F}_{n,\ell-1}}\stackrel{\mathrm{P}}{\to} a^2+2ab\Gamma+b^2\Sigma$, from the proof of Theorem \ref{thm:asymp}, the evaluation of $\Gamma$ below, and the local asymptotic normality of the model; the Lindeberg condition follows from $M_{n,\ell}(a,b)^2\leq 2a^2 D_{n,\ell}^2+2b^2\Delta_{n,\ell}^2$ and the Lindeberg conditions for $D_{n,\ell}$ and $\Delta_{n,\ell}$. Therefore, $\sum_{\ell=1}^n M_{n,\ell}(a,b)\inlaw\mathcal{N}(0,a^2+2ab\Gamma+b^2\Sigma)$ and the Cramér--Wold theorem entails that ${\bf T}_n\defin\sum_{\ell=1}^n \left(D_{n,\ell}, \Delta_{n,\ell} \right)^\top$ is asymptotically bivariate normal with mean zero and covariance matrix
$$
\begin{pmatrix}
1 & \Gamma \\
\Gamma & \Sigma
\end{pmatrix},
$$
where, due to the pairwise independence of the $\bX_{n,i}^\top\bX_{n,j}$'s and the orthogonality of the Gegenbauer polynomials,
\begin{align}
    \Gamma =&\;\lim_{n \to \infty} \frac{2 \sqrt{d_n}}{n^2} \sum_{\ell, \ell'=1}^n  \sum_{j=1}^{\ell-1}  \sum_{i=1}^{\ell'-1} {\rm Cov}(\psint(\bX_{n,i} ^\top \bX_{n,\ell'}), \bX_{n,j}^\top\bX_{n,\ell}) \nonumber \\
    =&\;\lim_{n \to \infty} \frac{2 \sqrt{d_n}}{n^2} \sum_{\ell, \ell'=1}^n  \sum_{j=1}^{\ell-1}  \sum_{i=1}^{\ell'-1} \E{\psint(\bX_{n,i} ^\top \bX_{n,\ell'}) \bX_{n,j}^\top\bX_{n,\ell}} \nonumber \\
    =&\; \lim_{n \to \infty} \frac{2 \sqrt{d_n}}{n^2} \sum_{\ell, \ell'=1}^n  \sum_{j=1}^{\ell-1}  \sum_{i=1}^{\ell'-1} \sum_{k=1}^\infty \lrp{1+\frac{2k}{d_n-1}}\tilde{v}^2_{k,d_n} \Ebig{C_k^{(d_n-1)/2}(\bX_{n,i} ^\top \bX_{n,\ell'}) \bX_{n,j}^\top\bX_{n,\ell}} \nonumber \\
    =&\;  \lim_{n \to \infty} \frac{ n(n-1) \lrp{(d_n-1)+{2}} \sqrt{d_n} \tilde{v}^2_{1,d_n}}{n^2} \E{(\bX_{n,1} ^\top \bX_{n,2})^2} \nonumber \\
    =&\;\lim_{n \to \infty} \frac{\sqrt{d_n} {v}^2_{1,d_n}}{\sqrt{2}\lrpbig{\sum_{k=1}^\infty v_{k,d_n}^4 d_{k,d_n}}^{1/2}}.\nonumber
\end{align}
The result then follows from a direct application of Le Cam's third Lemma \citep[see Theorem 6.6 in][]{vanderVaart1998}.
\end{proof}

\section{Required results}
\label{app:lemmas}

\begin{theorem}[Theorem 35.12 in \cite{Billingsley2012}] \label{thm:billingsley}
Let $D_{n,\ell}$, $\ell=1,\ldots,n$, $n=1,2,\ldots$, be a triangular array of random variables such that, for any $n$, $D_{n,1}, \ldots, D_{n,n}$ is a martingale difference sequence with respect to some filtration $\mathcal{F}_{n,1},\ldots, \mathcal{F}_{n,n}$. Assume that $\V{D_{n,\ell}}<\infty$ for any $(n,\ell)$. Letting $\sigma^2_{n,\ell}\defin\Ebig{D_{n,\ell}^{2} \mid \mathcal{F}_{n,\ell-1}}$ (with $\mathcal{F}_{n,0}$ being the trivial $\sigma$-algebra $\{\emptyset,\Omega\}$ for all $n$), further assume that
$$
\sum_{\ell=1}^{n} \sigma^2_{n,\ell} \stackrel{\mathrm{P}}{\to} 1
$$
as $n \to \infty$ and the Lindeberg condition
$$
\sum_{\ell=1}^{n} \E{D_{n,\ell}^2 1_{\{|D_{n,\ell}|>\varepsilon\}}} \to 0\quad\text{for all $\varepsilon>0$ as $n\to\infty$.}
$$
Then $\sum_{\ell=1}^{n} D_{n,\ell}\inlaw \mathcal{N}(0,1)$.
\end{theorem}

\begin{lemma} \label{lem:Ckm}
Let $\bX_{n,1}$ and $\bX_{n,2}$ be two independent random vectors uniformly distributed on $\Sdn$. For any integers $k,m\ge1$,
\begin{align*}
	\E{\lrb{C_{k}^{(d_n-1) / 2}(\bX_{n,1}^\top\bX_{n,2})}^{2m}}\asymp d_n^{mk}\quad\text{and}\quad d_{k,d_n}\sim \lrp{1+\frac{2k}{d_n-1}}\frac{d_n^k}{k!} \quad\text{as }d_n\to\infty,
\end{align*}
with $a_n\asymp b_n$ denoting that $a_n/b_n\to c\neq0$ as $n\to\infty$, and $a_n\sim b_n$ that $a_n/b_n\to 1$ as $n\to\infty$.
\end{lemma}

\begin{proof}[Proof of Lemma \ref{lem:Ckm}]
On the one hand, due to Equation 18.5.10 in \cite{NIST:DLMF},
\begin{align*}
	C_k^{(d_n-1)/2}(x)=\sum_{\ell=0}^{\lfloor k/2\rfloor} \frac{(-1)^\ell 2^{k-2\ell} \lrp{(d_n-1)/2}_{k-\ell}}{\ell!(k-2\ell)!} x^{k-2\ell}\indef \sum_{\ell=0}^{\lfloor k/2\rfloor} a_{k,\ell,d_n} x^{k-2\ell},
\end{align*}
where $(a)_k\defin\Gamma(a+k)/\Gamma(a)$ is the Pochhammer symbol. The coefficients $a_{k,\ell,d_n}$, $\ell=0,\ldots, \lfloor k/2\rfloor$, are such that
\begin{align}
	a_{k,\ell,d_n}=a_{k,\ell} \lrp{\frac{d_n-1}{2}}_{k-\ell}=a_{k,\ell} \frac{\Gamma\lrp{\frac{d_n-1}{2}+k-\ell}}{\Gamma\lrp{\frac{d_n-1}{2}}}\sim a_{k,\ell} \lrp{\frac{d_n-1}{2}+k-\ell}^{k-\ell}=O(d_n^{k-\ell}), \label{eq:O1}
\end{align}
where $a_{k,\ell}$ is a constant independent of $d_n$ and we have used Lemma B.3 in the Supplementary Material of \cite{Garcia-Portugues2023}: if $a\geq b$ and $\gamma>0$, then $\Gamma(\gamma d+a) / \Gamma(\gamma d+b) \sim(\gamma d+a)^{a-b}$ as $d \rightarrow \infty$.

On the other hand, from Lemma A.1 in \cite{Paindaveine2016}, we have that
\begin{align}
	\E{(\bX_{n,1}^\top\bX_{n,2})^{2m}}=\prod_{r=0}^{m-1} \frac{1+2r}{d_n+1+2r}\asymp d_n^{-m}. \label{eq:O2}
\end{align}

Joining \eqref{eq:O1} and \eqref{eq:O2} gives
\begin{align*}
	\E{\lrb{C_{k}^{(d_n-1) / 2}(\bX_{n,1}^\top\bX_{n,2})}^{2m}}=&\;\sum_{\ell_1,\ldots,\ell_{2m}=0}^{\lfloor k/2\rfloor} \lrb{\prod_{j=1}^{2m} a_{k,\ell_j,d_n}}\E{(\bX_{n,1}^\top\bX_{n,2})^{2mk-2\sum_{j=1}^{2m} \ell_j}}\\
	=&\;\sum_{\ell_1,\ldots,\ell_{2m}=0}^{\lfloor k/2\rfloor} O\lrp{d_n^{2mk-\sum_{j=1}^{2m} \ell_j}} O\lrp{d_n^{-(mk-\sum_{j=1}^{2m} \ell_j)}}\\
	=&\;O\lrp{d_n^{mk}}.
\end{align*}
This gives an upper bound for $\Ebig{\big\{C_{k}^{(d_n-1) / 2}(\bX_{n,1}^\top\bX_{n,2})\big\}^{2m}}$. For the lower bound, Jensen's inequality and \eqref{eq:ckd2} give $\Ebig{\big\{C_k^{(d_n-1)/2}(\bX_{n,1}^\top\bX_{n,2})\big\}^{2m}}\geq \Ebig{\big\{C_k^{(d_n-1)/2}(\bX_{n,1}^\top\bX_{n,2})\big\}^{2}}^m=((1+2k/(d_n-1))^{-2}d_{k,d_n})^m\asymp d_n^{mk}$, with the last equivalence following from \eqref{eq:dkdasymp} below.

Finally, \eqref{eq:dkd} and the previous asymptotic gamma expansion yield
\begin{align}
  d_{k,d_n}=\lrp{1+\frac{2k}{d_n-1}}\frac{\Gamma(d_n-1+k)}{\Gamma(d_n-1)k!}\sim\lrp{1+\frac{2k}{d_n-1}}\frac{d_n^k}{k!} \asymp d_n^k. \label{eq:dkdasymp}
\end{align}
\end{proof}

\fi

\ifsupplement

\newpage
\title{Supplement to ``High-dimensional Sobolev tests on hyperspheres''}
\setlength{\droptitle}{-1cm}
\predate{}%
\postdate{}%
\date{}

\author{Bruno Ebner$^{1}$, Eduardo Garc\'ia-Portugu\'es$^{2,4}$, and Thomas Verdebout$^{3}$}
\footnotetext[1]{Institute of Stochastics, Karlsruhe Institute of Technology (Germany).}
\footnotetext[2]{Department of Statistics, Universidad Carlos III de Madrid (Spain).}
\footnotetext[3]{Department of Mathematics, Université libre de Bruxelles (Belgium).}
\footnotetext[4]{Corresponding author. e-mail: \href{mailto:edgarcia@est-econ.uc3m.es}{edgarcia@est-econ.uc3m.es}.}
\maketitle

\begin{abstract}
These supplementary materials contain numerical experiments illustrating Theorems \ref{thm:asymp}--\ref{thm:LANstuff}.
\end{abstract}
\begin{flushleft}
	\small\textbf{Keywords:} Directional statistics; Normality; Rotational symmetry; Spherical symmetry; Uniformity.
\end{flushleft}

\appendix

\setcounter{section}{2}
\section{Numerical experiments}
\label{app:altern}

We conducted some numerical experiments to empirically evaluate the convergence of the test statistic $\sigma_n^{-1}T_n$ towards the asymptotic distributions stated in Theorems \ref{thm:asymp} and \ref{thm:LANstuff} when $n\to\infty$ and $d_n\to\infty$.

First, we considered $k_0$-Sobolev tests given by the sequence of weights $v_{k,d_n}=\delta_{k,k_0}$, for a fixed $k_0\geq1$. Under $\Hcal_{0,n}$, as given in \eqref{eq:asympCk}, we have
\begin{align}
    \tilde{T}_{n,k_0}\defin\frac{\sqrt{2\,k_0!}}{d_n^{k_0/2}n}\sum_{1\leq i<j\leq n} C_{k_0}^{(d_n-1)/2}(\bX_{n,i}^\top\bX_{n,j})\inlaw \mathcal{N}(0,1).\label{eq:asympCk2}
\end{align}
We simulated $M=10{,}000$ Monte Carlo samples $\bX_{n,1},\ldots, \bX_{n,n}$ uniformly distributed on $\Sdn$ for the pairs $(n,d_n)\in\{5,50,100,500\}^2$. Figures \ref{fig:C1}--\ref{fig:C2} show the histograms of sample realizations of $\tilde{T}_{n,k_0}$ for $k_0=1,3$, respectively, for varying $(n,d_n)$. The plots show that the histograms and the kernel density estimators (black curves) converge to the density of a $\mathcal{N}(0,1)$ (blue curve) as $\min(n,d_n)\to\infty$. The Kolmogorov--Smirnov (KS) test for a $\mathcal{N}(0,1)$ and the Lilliefors (L) test of normality are conducted to verify the convergence towards the limit distribution, and their $p$-values are reported for each scenario. For both figures, the $p$-values increasingly separate from zero as $\min(n,d_n)$ grows, considering that each $p$-value is asymptotically an independent uniform $[0,1]$ random variable.

\begin{figure}[h!]
    \centering
    \includegraphics[width=
    0.75\textwidth]{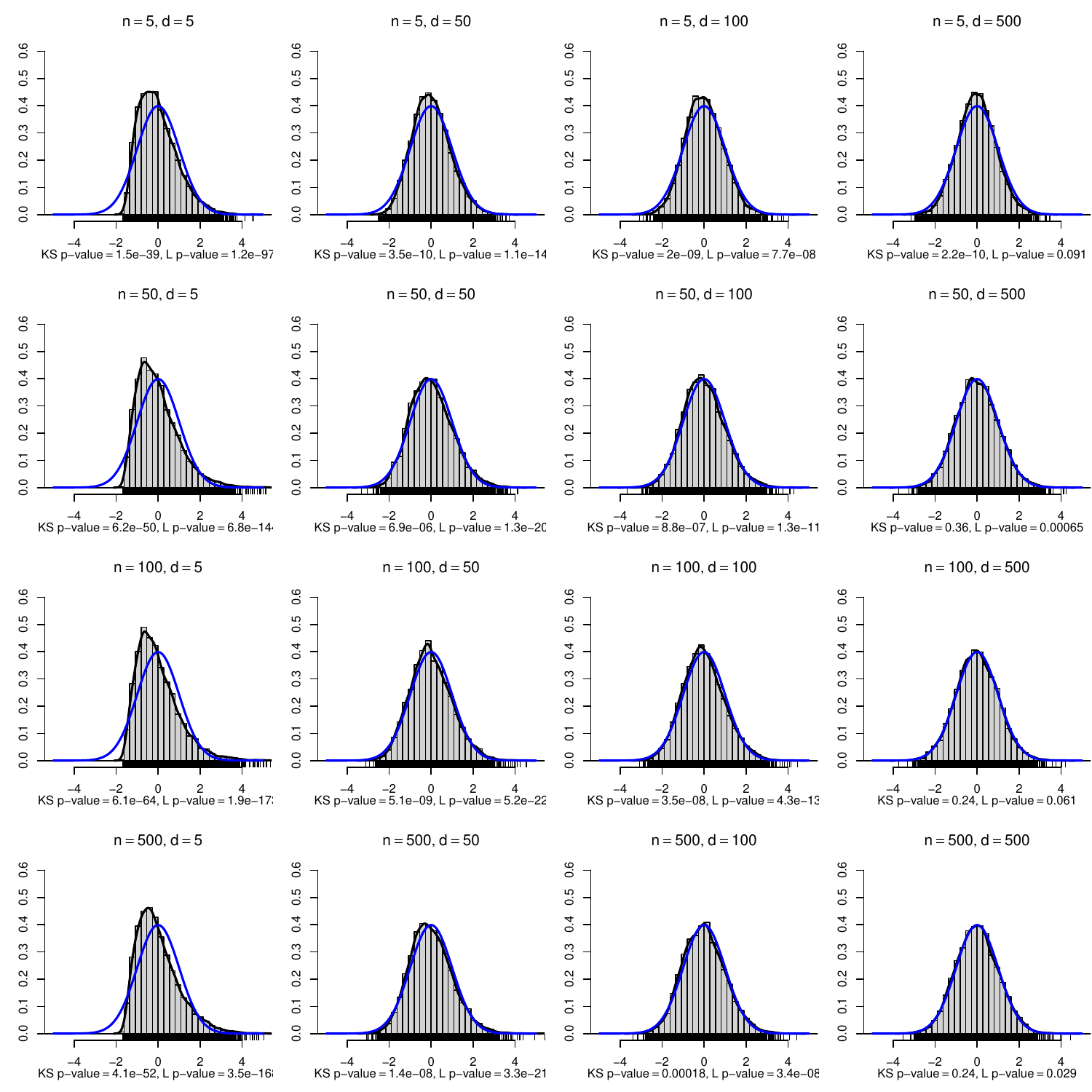}
    \caption{\small Empirical evaluation of $\tilde{T}_{n,k_0}\inlaw \mathcal{N}(0,1)$ under $\Hcal_{0,n}$ for $k_0=1$.}
    \label{fig:C1}
\end{figure}

Second, we repeated the same experiment for $k_0$-Sobolev tests, but now with samples $\bX_{n,1},\ldots,\bX_{n,n}$ simulated under $\Hcal_{1,n}$ defined by an integrated von Mises--Fisher distribution \eqref{eq:ivMF} on $\Sdn$ with concentration $\kappa_n = \tau d_n^{3/4}/\sqrt{n}$ for $\tau>0$. According to Theorem \ref{thm:LANstuff} and subsequent comments, \eqref{eq:asympCk2} becomes
$$
\tilde{T}_{n,k_0}\inlaw\mathcal{N}(\Gamma \tau^2,1)
$$
under $\Hcal_{1,n}$, with $\Gamma=(1/\sqrt{2}) 1_{\{k_0=1\}}$. Setting $\tau^2=\sqrt{2}$ yields the asymptotic non-null distribution $\mathcal{N}(1_{\{k_0=1\}},1)$. Figure \ref{fig:rem1} shows the convergence of the histogram and kernel density estimator of the sample realizations of $\tilde{T}_{n,k_0}$ for $k_0=1$ toward the density of a $\mathcal{N}(1_{\{k_0=1\}},1)$ (red curve). The Kolmogorov--Smirnov tests are now conducted for this distribution. For $k_0=1$, where asymptotic power is present, convergence toward the limit distribution is slower than in the null hypothesis, particularly in terms of the variance.

Third, we considered finite Sobolev tests with $k_0>1$ positive weights for the same integrated von Mises--Fisher alternatives. According to the comments after Theorem \ref{thm:LANstuff}:
\begin{enumerate}
    \item $\Gamma=0$ for the sequence of weights $v_{k,d_n}=1_{\{k\leq k_0\}}$.
    \item $\Gamma=1/\sqrt{2}$ for the sequence of weights $v_{k,d_n}=\delta_{k,1}+[k! d_n^{-k}]^{1/4}1_{\{1<k\leq k_0\}}$.
\end{enumerate}
Considering $\tau^2=\sqrt{2}$, then $\sigma_n^{-1}T_n\inlaw \mathcal{N}(0,1)$ for the first case and $\sigma_n^{-1}T_n\inlaw \mathcal{N}(1,1)$ for the second. Figures \ref{fig:rem3} and \ref{fig:rem4} illustrate these behaviors when $k_0=3$ (i.e., only the first weights are positive).

\begin{figure}[!htb]
    \centering
    \includegraphics[width=0.75\textwidth]{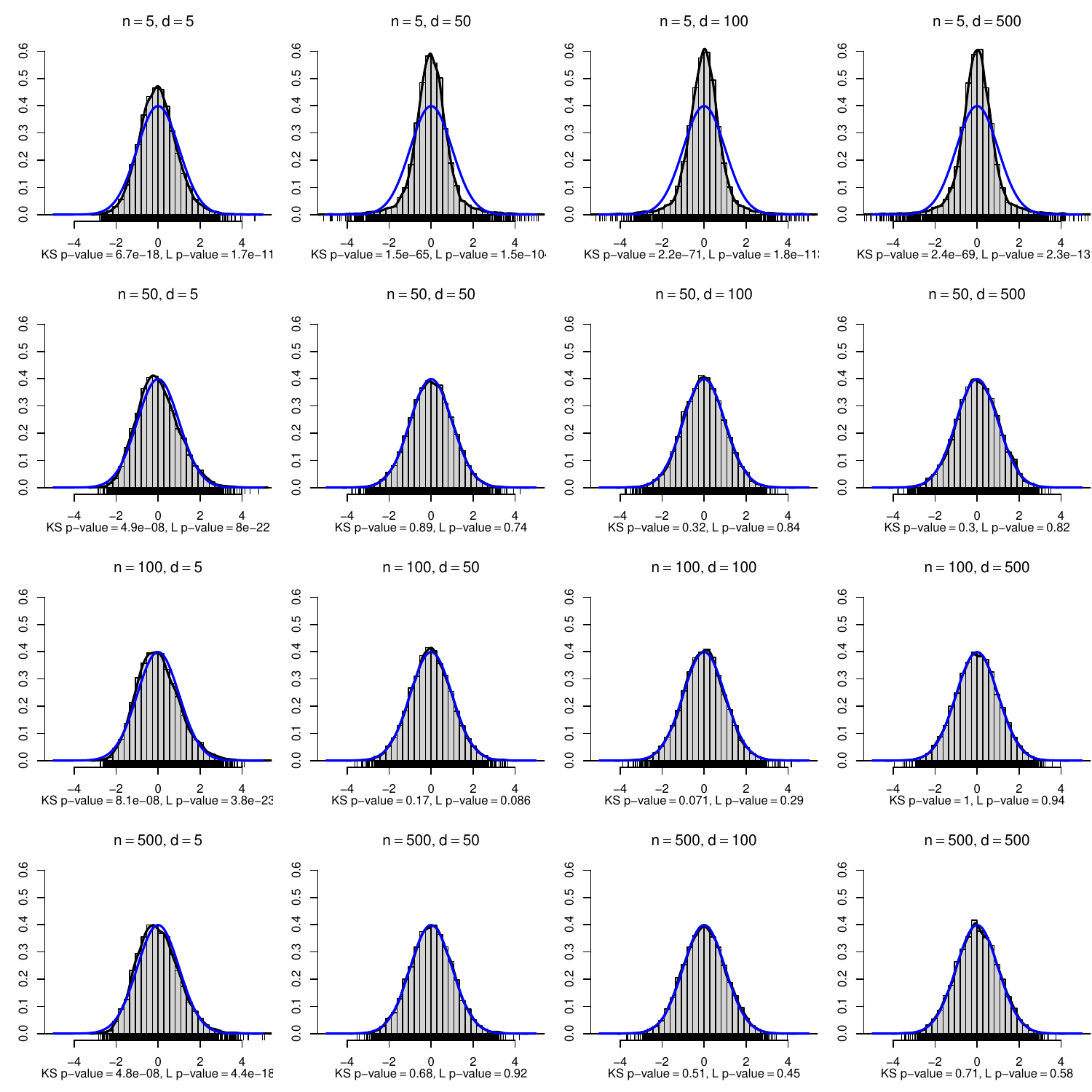}
    \caption{\small Empirical evaluation of $\tilde{T}_{n,k_0}\inlaw \mathcal{N}(0,1)$ under $\Hcal_{0,n}$ for $k_0=3$.}
    \label{fig:C2}
\end{figure}

\begin{figure}[!htb]
    \centering
    \includegraphics[width=0.75\textwidth]{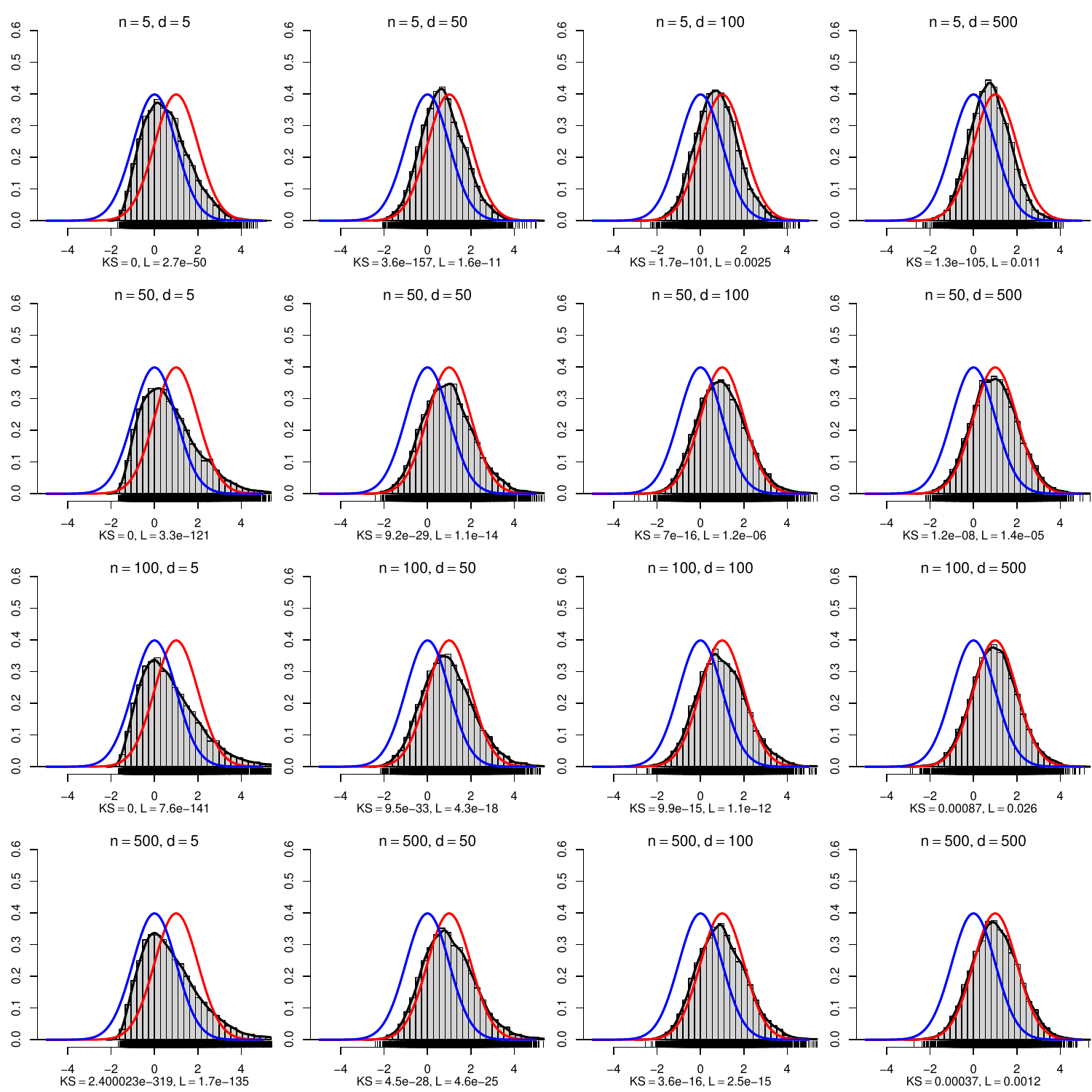}
    \caption{\small Empirical evaluation of $\tilde{T}_{n,k_0}\inlaw \mathcal{N}(1,1)$ under $\Hcal_{1,n}$ for $k_0=1$.}
    \label{fig:rem1}
\end{figure}

\begin{figure}[!htb]
    \centering
    \includegraphics[width=0.75\textwidth]{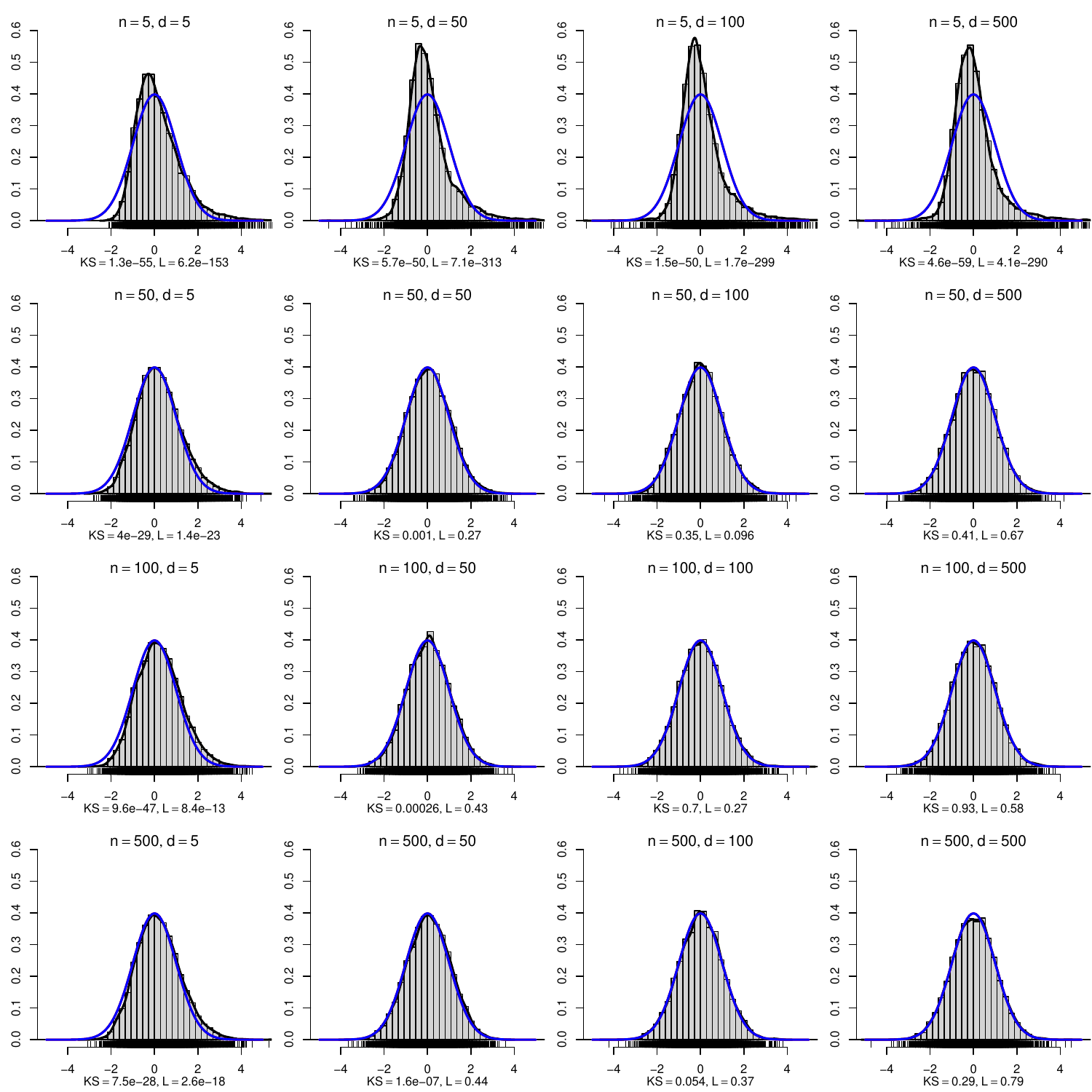}
    \caption{\small Empirical evaluation of $\sigma_n^{-1}T_n\inlaw \mathcal{N}(0,1)$ under $\Hcal_{1,n}$ with $v_{k,d_n}=1_{\{k\leq k_0\}}$ and $k_0=3$.}
    \label{fig:rem3}
\end{figure}

\begin{figure}[!htb]
    \centering
    \includegraphics[width=0.75\textwidth]{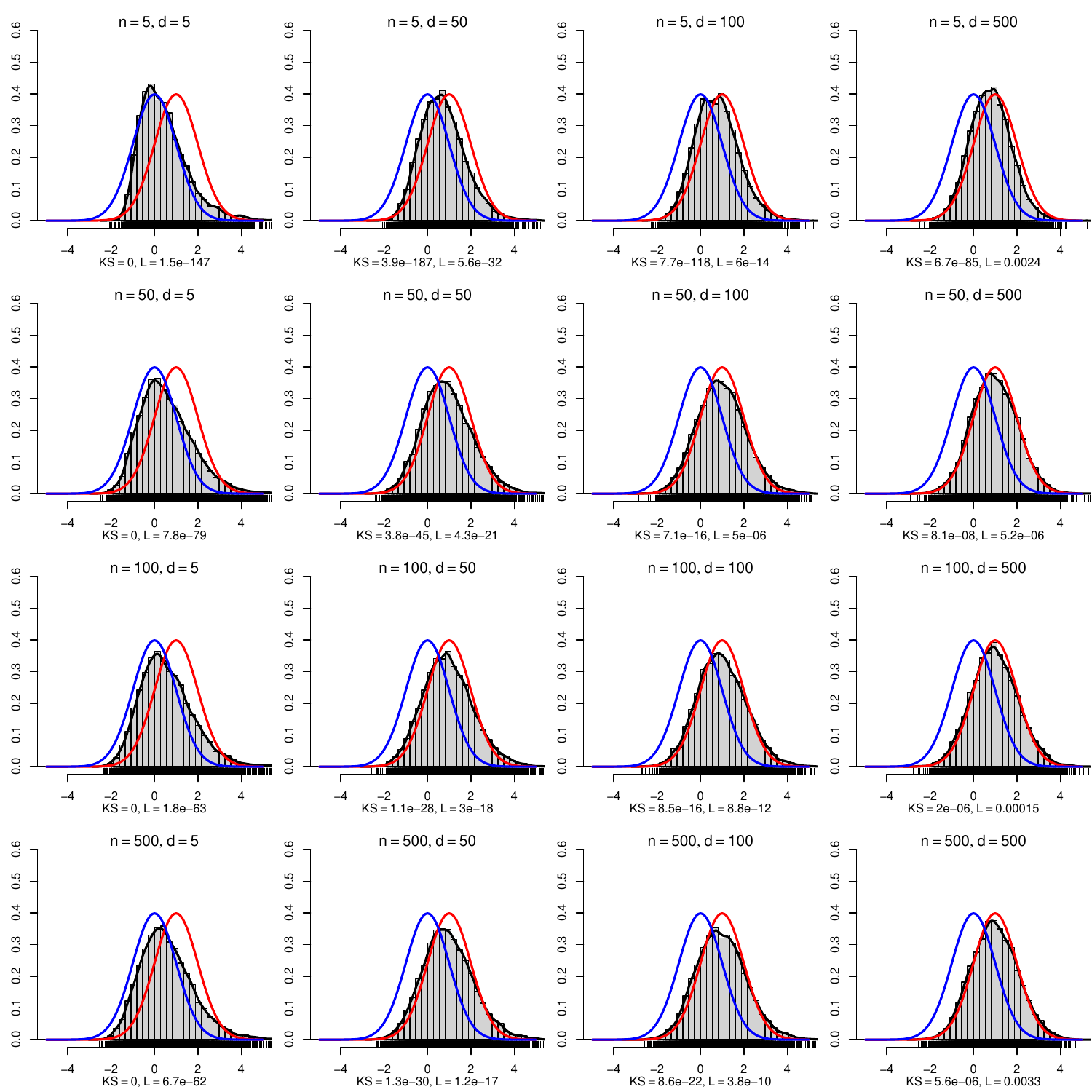}
    \caption{\small Empirical evaluation of $\sigma_n^{-1}T_n\inlaw \mathcal{N}(1,1)$ under $\Hcal_{1,n}$ with $v_{k,d_n}=\delta_{k,1}+[k! d_n^{-k}]^{1/4}1_{\{1<k\leq k_0\}}$ and $k_0=3$.\!\!\!}
    \label{fig:rem4}
\end{figure}

\fi

\end{document}

%% file: preamble_arXiv.tex
\usepackage[utf8]{inputenc}
\usepackage[T1]{fontenc}

\usepackage{amsthm}
\usepackage{amsmath}
\usepackage{amsfonts}
\usepackage{amssymb}

\usepackage{graphicx} 
\usepackage{graphbox}
\usepackage{float}
\usepackage{booktabs}
\usepackage{multirow}
\usepackage{rotating}
\usepackage{array}
\usepackage{xcolor}
\usepackage{subcaption}
\usepackage[ruled]{algorithm2e}

\newcolumntype{C}[1]{>{\centering\arraybackslash}p{#1}}
\newcolumntype{R}[1]{>{\raggedleft\let\newline\\\arraybackslash\hspace{0pt}}m{#1}}
\newcolumntype{L}[1]{>{\raggedright\let\newline\\\arraybackslash\hspace{0pt}}m{#1}}

\usepackage{breakcites}

\usepackage{nowidow}

\usepackage{enumitem}

\usepackage[top=1.5cm,bottom=2cm,right=2.5cm,left=2.5cm]{geometry}
\usepackage{titling}

\usepackage{natbib}

\usepackage{ifplatform}

\usepackage{textcomp}
\usepackage{bbm}

\usepackage{comment}

\ifwindows
\fi

\allowdisplaybreaks

\newcommand{\lc}{\left[}
\newcommand{\rc}{\right]}

\newcommand{\R}{\mathbb{R}}
\newcommand{\Sd}{\mathbb{S}^d}
\newcommand{\Sdn}{\mathbb{S}^{d_n}}
\newcommand{\Sdnm}{\mathbb{S}^{d_n-1}}

\newcommand{\rd}{\mathrm{d}}

\newcommand{\bx}{\boldsymbol{x}}

\newcommand{\bX}{\boldsymbol{X}}
\newcommand{\bY}{\boldsymbol{Y}}

\newcommand{\bO}{\boldsymbol{O}}
\newcommand{\bmu}{\boldsymbol\mu}
\newcommand{\bu}{\boldsymbol{u}}
\newcommand{\bv}{\boldsymbol{v}}

\newcommand{\bt}{\boldsymbol{t}}

\newcommand{\zero}{\boldsymbol{0}}
\newcommand{\bxi}{\boldsymbol\xi}

\newcommand{\btheta}{\boldsymbol\theta}

\newcommand{\brho}{\boldsymbol\rho}

\newcommand{\bOmega}{\boldsymbol\Omega}

\newcommand{\bSigma}{\boldsymbol\Sigma}

\newcommand{\bGamma}{\boldsymbol\Gamma}

\newcommand{\Hcal}{\mathcal{H}}

\newcommand{\bB}{\boldsymbol{B}}

\newcommand{\bI}{\boldsymbol{I}}
\newcommand{\bS}{\boldsymbol{S}}

\newcommand{\inlaw}{\rightsquigarrow}

\newcommand{\psin}{\psi_n}
\newcommand{\psint}{\tilde{\psi}_n}
\newcommand{\defin}{:=}
\newcommand{\indef}{=:}

\newcommand{\lrp}[1]{\left(#1\right)}
\newcommand{\lrpbig}[1]{\big(#1\big)}

\newcommand{\lrb}[1]{\left\{#1\right\}}
\newcommand{\Prob}[1]{\mathrm{P}\lc #1\rc}
\newcommand{\E}[1]{\mathrm{E}\lc #1\rc}
\newcommand{\Ebig}[1]{\mathrm{E}\big[ #1\big]}
\newcommand{\EBig}[1]{\mathrm{E}\Big[ #1\Big]}
\newcommand{\V}[1]{\mathrm{Var}\lc #1\rc}
\newcommand{\Vbig}[1]{\mathrm{Var}\big[ #1\big]}

\newcommand{\om}[1]{\omega_{#1}}

\DeclareFontFamily{OT1}{pzc}{}
\DeclareFontShape{OT1}{pzc}{m}{it}{<-> s * [1.10] pzcmi7t}{}
\DeclareMathAlphabet{\mathpzc}{OT1}{pzc}{m}{it}

\theoremstyle{plain}

\newtheorem{theorem}{Theorem}[section]

\newtheorem{corollary}{Corollary}[section]
\newtheorem{lemma}[theorem]{Lemma}
\theoremstyle{definition}

\newtheorem{remark}[theorem]{Remark}

\newcommand{\pct}{\%}